\documentclass[10pt]{amsart} \usepackage{amssymb, verbatim} 

\newcommand{\Mbar}{\overline{M}}
\def\R{{\hbox{\bf R}}}

\def\Re{{\hbox{Re}}}

\def\R{{\hbox{\bf R}}}
\def\C{{\hbox{\bf C}}}

\def\Z{{\hbox{\bf Z}}}

\def\eps{\varepsilon}

\def \endprf{\hfill  {\vrule height6pt width6pt depth0pt}\medskip}
\def\emph#1{{\it #1}}
\def\textbf#1{{\bf #1}}
\def\divider#1{$\bullet\quad${\bf #1}}
\def\sc{\operatorname{sc}}
\def\Psisc{\Psi_{\sc}} 
\def\Hess{{\operatorname{Hess}}}

\def\f12{\frac{1}{2}}

\DeclareMathOperator{\Op}{{Op}}

\newcommand{\bea}{\begin{eqnarray}}
\newcommand{\eea}{\end{eqnarray}}
\newcommand{\beaa}{\begin{eqnarray*}}
\newcommand{\eeaa}{\end{eqnarray*}}

\def\divider#1{$\bullet\quad${\bf #1}}
 
\theoremstyle{plain}
  \newtheorem{theorem}[subsection]{Theorem}
  
  \newtheorem{proposition}[subsection]{Proposition}
  \newtheorem{lemma}[subsection]{Lemma}

\theoremstyle{remark}

\theoremstyle{definition}

\include{psfig}

\begin{document}

\title[Long-time decay for Schr\"odinger]{Long-time decay estimates 
for the Schr\"odinger equation on manifolds} 

\author{Igor Rodnianski}
\address{Department of Mathematics, Princeton University, Princeton 
NJ 08544}
\email{ irod@math.princeton.edu}

\author{Terence Tao}
\address{Department of Mathematics, UCLA, Los Angeles CA 90095-1555}
\email{ tao@math.ucla.edu}


\vspace{-0.3in}

\begin{abstract}  In this paper we develop a quantitative version of Enss' method to establish global-in-time decay estimates for
solutions to Schr\"odinger equations on manifolds.  To simplify the exposition we shall only consider
Hamiltonians of the form $H := - \frac{1}{2} \Delta_M$,
where $\Delta_M$ is the Laplace-Beltrami operator on a manifold $M$ which is a smooth compact perturbation of three-dimensional
Euclidean space $\R^3$ which obeys the non-trapping condition.  We establish a global-in-time local 
smoothing estimate for the Schr\"odinger equation $u_t = -iHu$.  The main novelty here is
the global-in-time aspect of the estimates, which forces a more detailed analysis on
the low and medium frequencies of the evolution than in the local-in-time theory.
In particular, to handle the medium frequencies we require the RAGE theorem (which reflects the fact
that $H$ has no embedded eigenvalues), together with a quantitative version of Enss' method decomposing the solution asymptotically
into incoming and outgoing components, while to handle the low frequencies we need a Poincar\'e-type inequality 
(which reflects the fact that $H$ has no eigenfunctions or resonances at zero).
\end{abstract}

\maketitle

\section{Introduction}

Let $(M,g) = (\R^3,g)$ be a compact perturbation of Euclidean space\footnote{The analysis we give here also extends to higher dimensions $n>3$, and are in fact slightly easier in those cases, but for simplicity of exposition we restrict our attention to the physically important three-dimensional case
and to compact perturbations.  The hypothesis that $M$ is topologically $\R^3$ is technical but in any event is forced upon us by the non-trapping hypothesis which we introduce below; see \cite{Thorbergsson}.} $\R^3$, thus $M$ is $\R^3$ endowed with a smooth metric
$g$ which equals the Euclidean metric outside of a Euclidean ball $B(0,R_0) := \{ x \in \R^3: |x| \leq R_0\}$ for some fixed $R_0$.  We consider 
smooth solutions to the Schr\"odinger equation
\begin{equation}\label{schrodinger}
u_t = -iHu
\end{equation}
where for each time $t$, $u(t): M \to \C$ is a Schwartz function on $M$, and $H$ is the Hamiltonian operator
$$ H := - \frac{1}{2} \Delta_M$$
where $\Delta_M := \nabla^j \nabla_j$ is the Laplace-Beltrami operator (with $\nabla^j$ denoting covariant 
differentiation with respect to the Levi-Civita connection, in contrast with the Euclidean partial derivatives $\partial_j$).
Note that $H$ is positive definite and self-adjoint with respect to the natural inner product
\begin{equation}\label{inner-def}
 \langle u, v \rangle_{L^2(M)} := \int_M u(x) \overline{v(x)}\ dg(x),
\end{equation}
where $dg :=\sqrt{\det g_{ij}(x)} dx$ is the standard volume element induced by the metric $g$.  
In fact the spectrum of $H$ consists entirely of absolutely continuous spectrum on the positive real axis $[0,+\infty)$,
in particular $H$ has no eigenvalues or resonances at any energy.  In particular $H$ enjoys a standard functional calculus
on $L^2(M)$, the Hilbert space associated to the inner product, and one can define the homogeneous Sobolev norms 
$\| u \|_{\dot H^s(M)} := \| H^{s/2} u \|_{L^2(M)}$ on $M$ using fractional powers of $H$ for all $-1 \leq s \leq 1$.  (One can of course define
these norms for other $s$ also, but there are some technicalities when $s$ is too negative that we do not wish to address here).  Thus for instance
$\| u \|_{\dot H^0(M)} = \| u \|_{L^2(M)}$ and
\begin{equation}\label{h1}
\| u \|_{\dot H^1(M)}^2 = \langle u, Hu \rangle_{L^2(M)} = \frac{1}{2} \int_M |\nabla u|^2_g\ dg.
\end{equation}

It is well known (see \cite{cks}) that for any time $t_0$ and any Schwartz initial data $u_0$, there exists a unique global-in-time
Schwartz solution $u: \R \times M \to \C$ to \eqref{schrodinger} with initial data $u(0) = u_0$, indeed
we have $u(t) = e^{-it H} u_0$.  In this paper we develop a quantitative variant of Enss' method to obtain a new \emph{global-in-time}
local smoothing estimate for such solutions; to avoid needless technicalities we shall always
restrict ourselves to Schwartz solutions.  The methods here extend to more general classes of Hamiltonians than those considered here,
for instance they can handle asymptotically flat manifolds of dimension $n \geq 3$ as well as short-range potentials provided that there are no resonances or eigenfunctions
at zero; however to simplify the exposition we have chosen to restrict attention to the simple case of zero potential and compact perturbations
of three-dimensional Euclidean space.  We shall pursue more general Hamiltonians in \cite{rt:limiting} using a somewhat different method (based on limiting absorption
principles).

Let us begin by recalling some earlier results.  In the case where $M$ is Euclidean space $\R^n$, so that
$H = H_0 := - \frac{1}{2} \Delta_{\R^3}$ is just the free Hamiltonian, then we have the 
well-known global-in-time local smoothing estimate (see \cite{constantin}, \cite{sjolin}, \cite{vega}),
which we shall phrase as
\begin{equation}\label{eq:compact}
 \int_{-\infty}^{\infty} \| \langle x \rangle^{-1/2-\sigma} \nabla e^{-itH_0} u_0 \|_{L^2(\R^3)}^2
+ \| \langle x \rangle^{-3/2-\sigma} e^{-itH_0} u_0 \|_{L^2(\R^3)}^2\ dt
 \leq C_{\sigma} \| u_0 \|_{\dot H^{1/2}(\R^3)}^2
\end{equation}
for any $\sigma > 0$ and any Schwartz initial data $f$, where $\langle x \rangle := (1 + |x|^2)^{1/2}$.
It is well known that the condition on $\sigma$ is sharp (there is a logarithmic divergence in the left-hand side
when $\sigma = 0$).

Now we suppose that $(M,g) = (\R^3,g)$ is a compact perturbation of Euclidean space obeying the non-trapping condition.  
It is known (see \cite{cks}, \cite{doi}) that one
has the \emph{local-in-time} local smoothing estimate
\begin{equation}\label{ls-manifold}
 \int_I \| \langle x \rangle^{-1/2-\sigma} \nabla e^{-it H} u_0 \|_{L^2(M)}^2
+ \| \langle x \rangle^{-3/2-\sigma} e^{-it H} u_0 \|_{L^2(M)}^2\ dt
 \leq C_{\sigma,I,M} \| u_0 \|_{\dot H^{1/2}(M)}^2
\end{equation}
for all \emph{compact} time intervals $I \subset \R$, where the constant on the right-hand side is allowed to
depend on $I$, if and only if the manifold $M$ is \emph{non-trapping}, i.e. every geodesic $s \mapsto x(s)$ in $M$ eventually goes to spatial infinity as have a theory for the Schr\"odinger equation as $s \to \pm \infty$.  To understand why the non-trapping condition is necessary, observe that the localization in time means that the low and medium energies are easily controlled (for instance by using the Sobolev embedding $\dot H^{1/2}(M) \subseteq L^3(M)$ and H\"older's inequality, noting that at low and medium energies $\nabla e^{-itH} u_0$ is also in $\dot H^{1/2}(M)$), and so one only needs to understand the evolution of the high energies, which evolve semi-classically.  The semi-classical limit of the estimate \eqref{ls-manifold} is the 
estimate $\int_\R \langle x(s) \rangle^{-1-2\sigma} |\dot x(s)|^2\ ds \leq C_{\sigma,M} |\dot x(0)|$ for any geodesic $s \mapsto x(s)$,
which is easily seen to hold if and only if the manifold is non-trapping (this also explains the requirement that $\sigma > 0$).  
Of course, one needs semi-classical tools such as pseudo-differential operators and the positive commutator method 
in order to make this argument rigorous; see \cite{cks}, \cite{doi} (or Section \ref{highfreq-sec}) for more details.

The main result of this paper is to unify the global-in-time Euclidean estimate \eqref{eq:compact} with the local-in-time manifold estimate
\eqref{ls-manifold} as follows.

\begin{theorem}\label{main}  Let $M$ be a smooth compact perturbation of $\R^3$ which is non-trapping and which is smoothly diffeomorphic to $\R^3$.  Then for any Schwartz solution $u(t,x)$ to \eqref{schrodinger} and any $\sigma > 0$ we have
\begin{equation}\label{ls-manifold-global}
 \int_\R \| \langle x \rangle^{-1/2-\sigma} \nabla e^{-it H} u_0 \|_{L^2(M)}^2
+ \| \langle x \rangle^{-3/2-\sigma} e^{-it H} u_0 \|_{L^2(M)}^2\ dt
 \leq C_{\sigma,M} \| u_0 \|_{\dot H^{1/2}(M)}^2.
\end{equation}
In other words, the constant $C_{\sigma,I,M}$ in
\eqref{ls-manifold} can be taken to be independent of the interval $I$.
\end{theorem}

As mentioned earlier, this is not the strongest result that one could obtain with this method; our purpose here is merely to illustrate
a model example in which the method applies.  By using different methods (in particular the limiting absorption principle)
we were able to obtain results for more general Hamiltonians, see \cite{rt:limiting}.  However we believe the method we present here still
has some merit, in particular it is conceptually straightforward and seems to have some hope of generalizing to more ``time-dependent''
or ``non-linear'' situations in which spectral theory tools are less useful.  For instance, the methods here were already used in \cite{tao:focusing}
to obtain new results on the asymptotic behavior of focusing non-linear Schr\"odinger equations, at least in the spherically symmetric case.

We now informally discuss the proof of the theorem.  First note that this estimate is already proven when $M$ is Euclidean space $\R^3$,
and we are only considering manifolds $M$ which are compact perturbations of Euclidean space.  Thus it is reasonable to expect that the only
difficulties in proving \eqref{ls-manifold-global} will arise from the compact region $B(0,R_0)$,
and indeed by modifying the proof of \eqref{eq:compact} we will be able to show that a global-in-time 
local smoothing estimate in a slight enlargement $B(0,R)$ of $B(0,R_0)$ will automatically imply the full estimate \eqref{ls-manifold-global}.  Thus we may (heuristically, at least) restrict our attention to a ball such as $B(0,R)$.

Next, we use the spectral theorem to decompose the evolution \eqref{schrodinger} into low-energy, medium-energy, and high-energy components.  The high energy components turn out to be treatable by the same positive commutator arguments used to prove the local-in-time estimate \eqref{ls-manifold}.  From
a heuristic viewpoint, this is because high-energy components propagate very quickly, and thus only linger in the compact region $B(0,R)$ for a very
short period of time, after which they escape into the Euclidean region of $M$ and never return to $B(0,R)$ again.  Of course the non-trapping hypothesis is essential here.  At a more technical level, the reason why we can adapt the arguments used to prove \eqref{ls-manifold} is that the error terms generated by the positive commutator method are lower order than the main term, and can thus be absorbed by the main term in the high-energy regime
even when the time interval $I$ is unbounded.

The low-energy components are easy to treat, but for a different reason, namely that there is an ``uncertainty principle'' that shows that solutions which have extremely low energy cannot be concentrated entirely in the compact region $K$, and thus the low energies cannot be the dominant component
to the local smoothing estimate inside this region $K$.  This is ultimately a reflection of the well-known fact that a Hamiltonian $H = -\frac{1}{2} \Delta_M$ with no potential does not have any resonances or bound states at zero or negative energies.

The most interesting component to treat is the medium energies regime, which requires new methods.  These are energies which are not high enough to
behave semi-classically and escape the compact region $K$ by means of the non-trapping hypothesis, but which are not low enough to escape the compact region $K$ by means of the uncertainty principle.  Fortunately there is a third mechanism by which we can force solutions of the Schr\"odinger equation to escape the compact region $K$, namely the Ruelle-Amrein-Gorgescu-Enss (RAGE) theorem.  The point is that (as is well known) $H$ contains no embedded eigenfunctions 
in the medium energy portion of the spectrum (or indeed anywhere in the spectrum), and so the RAGE theorem then ensures that any given 
solution to \eqref{schrodinger} must eventually
vacate the region $K$ after some time $T$.  
We recall the abstract version of the
RAGE Theorem (see e.g. \cite{reed}).

\begin{theorem}[RAGE]\label{rage-thm}
Let $H$ be a self-adjoint operator on a Hilbert space ${\mathcal H}$. If 
$C:{\mathcal H}\to {\mathcal H}$ is a compact operator and 
$u_0$ lies in the continuous subspace of $H$, then
$$
\lim_{T\to \infty} \frac 1T \int_0^T \|Ce^{-itH} u_0\|_H^2=0
$$
\end{theorem}

The RAGE theorem played a crucial role in Enss' approach to scattering for 
Schr\"odinger operators $H=H_0+V$, \cite{enss}. In fact our treatment of medium energies 
can be viewed as a quantitative version of the Enss' method in a sense that we use 
both the RAGE type inequality and the decomposition into incoming and outgoing
waves, also a major part of the Enss' work, but derive an a priori estimate, as opposed
to a qualitative result about completeness of the wave operators. The other major 
difference is that in scattering theory, density arguments allow one to consider
{\it only} medium energies and compactly supported data, whereas we must necessarily treat all 
energy ranges and allow our data to have arbitrary support.

One might object that the RAGE theorem is ``qualitative'' in nature, in that the time $T$ required for a solution to leave $K$ depends on the choice of solution and thus need not be uniform.  However because we have localized the solution in both frequency (to medium energies) and position (to the region $K$), and because we can use linearity to normalize the $H^{1/2}(M)$ norm of $u$, the solution is in fact effectively contained in a compact region of (quantum) phase space.  Because of this, one can make the time $T$ required for a solution to
leave $K$ to be uniform for all medium-energy solutions $u$.

There is however still a remaining difficulty for medium energy solutions, which is that once a solution leaves the compact region $K$ one needs to
ensure that it does not return back to $K$, since a solution which periodically left and then returned to $K$ would eventually contribute an
infinite left-hand side to \eqref{ls-manifold-global}.  To resolve this we introduce a quantitative version of Enss' method.  The starting point is the observation that, in the exterior of the domain $K$, any function can be decomposed into ``outgoing'' and ``incoming'' components (very roughly speaking, this corresponds to the spectral projections $\chi_{[0,+\infty)}(i \partial_r)$ and $\chi_{(-\infty,0]}(i \partial r)$ where $\partial_r$ is
the radial derivative, although we shall not perform these projections directly due to the singular behavior of the operator $i \partial_r$).
Outgoing components will evolve towards spatial infinity as $t \to +\infty$, whereas incoming components will evolve towards spatial infinity
as $t \to -\infty$.  In particular in both cases the solution will not encounter the compact region $K$ and the evolution is essentially Euclidean in nature.  

Now suppose that a medium energy solution to \eqref{schrodinger} is localized to $K$ at some time $t_0$.  By the quantitative RAGE theorem, at some later time $t_0 + T$ the solution has mostly vacated the region $K$.  By Enss' decomposition it can have either outgoing or incoming components.  But one can show that there is almost no incoming component, because if we evolved backwards in time from $t_0+T$ back to $t_0$ we see that the incoming component would
have evolved back to a region far away from $K$, and thus orthogonal to the initial data.  Thus at time $t_0 + T$, the solution consists almost primarily of outgoing components.  But then by the preceding discussion this means that the solution will continue to radiate to spatial infinity
for times after $t_0$, and thus never return to $K$.  To summarize, we have shown that any component of a medium energy solution which is located in $K$ at time $t_0$, will eventually radiate to spatial infinity as $t \to +\infty$; a similar argument also handles the $t \to -\infty$ evolution.  Combining this with the finite energy of $u$ (note that in the medium energy regime all Sobolev norms are equivalent), one can obtain
the estimate \eqref{ls-manifold-global}.

The global-in-time local smoothing estimate in Theorem \ref{main} has a number of consequences; for instance by combining it
with the arguments of Staffilani and Tataru \cite{st} one can obtain global-in-time Strichartz estimates for compact non-trapping perturbations
of Euclidean space, which then can be used to transfer some local and global existence results for non-linear Schr\"odinger equations on
Euclidean space, to the setting of compact non-trapping perturbations.  We will not discuss these (fairly standard) 
generalizations here, but see \cite{rt:limiting} for further discussion.

I.R. is a Clay Prize Fellow and supported in part by the NSF grant DMS-01007791. T.T. is supported in part by a grant from the Packard Foundation.
We thank Nicolas Burq for helpful comments.

\section{Notation}

Throughout this paper, the manifold $M$, the radius $R_0$, and the
exponent $\sigma$ will be fixed. All constants $C$ will be allowed to depend on $\sigma$ and $M$ (and hence on $R_0$).  If $C$ needs to depend on other parameters we will indicate this by subscripts.

Suppose $A_0, A_1, \ldots, A_k$ are real parameters.  We use $O_{A_1,\ldots,A_k}(1)$ to denote
any quantity depending on $A_1,\ldots,A_k$ and possibly some other quantities, which is bounded in magnitude
by some constant $C_{A_1,\ldots,A_k}$.  For some fixed $c$ (usually $c=0$ or $c=\infty$), we 
also use $o_{A_0 \to c; A_1,\ldots,A_k}(1)$ to denote a quantity
depending on $A_0, A_1,\ldots, A_k$ and possibly some other parameters, which is bounded in magnitude by some quantity
$F_{A_1,\ldots,A_k}(A_0)$ such that 
$$ \lim_{A_0 \to c} F_{A_1,\ldots,A_k} = 0 \hbox{ for all choices of } A_1,\ldots, A_k.$$
Often $k$ will be zero, in which case the above notations simply read $O(1)$ and $o_{A_0 \to c}(1)$ respectively.
We also abbreviate $O_{A_1,\ldots,A_k}(1) X$ and $o_{A_0 \to c; A_1,\ldots,A_k}(1) X$
as $O_{A_1,\ldots,A_k}(X)$ and $o_{A_0 \to c; A_1,\ldots,A_k}(X)$ respectively.  

In the absence of parentheses, we read operators from right to left.  Thus for instance, $\nabla f g$ denotes
the function $\nabla(fg)$ rather than $(\nabla f)g$.

We use the usual summation conventions on indices, and use $g^{jk}$ to denote the dual metric to $g_{jk}$ on
the cotangent bundle, and use $g$ to raise and lower indices in the usual manner.
We use $\nabla_j$, $\nabla^j$ to denote the usual covariant derivatives with respect to the Levi-Civita connection
on $M$; these can be applied to any tensor field; we use $|\nabla f|_g = \sqrt{\nabla^j f \overline{\nabla_j f}}$ to denote the magnitude of the gradient with respect to the metric $g$, and $|\nabla f|$ to denote the Euclidean magnitude of the gradient.  
Since $\nabla g = 0$, the indices of covariant derivatives
can be raised and lowered freely, thus for instance $\Delta_M = \nabla^j \nabla_j = \nabla_j \nabla^j$.  Also
these covariant derivatives are anti-selfadjoint with respect to the inner product \eqref{inner-def}, and
thus we can integrate by parts using these derivatives freely.  

As the operator $H$ is self-adjoint and has spectrum on $[0,\infty)$ we can construct spectral multipliers
$f(H)$ for any measurable function $f: [0,\infty) \to \C$ of at most polynomial growth; in particular
we can define fractional powers $H^{s/2}$ and $(1+H)^{s/2}$, as well as Schr\"odinger propagators $e^{-itH}$
and Littlewood-Paley type operators on $H$.  These spectral multipliers commute with each other, and are bounded
on $L^2$ if their symbol $f$ is bounded.  

\section{Overview of proof}

We now begin the proof of Theorem \ref{main}.  The first step will be to show that one can freely
pass back and forth between the slowly decaying weight
$\langle x \rangle^{-1/2-\sigma}$ in \eqref{ls-manifold-global} and a suitably chosen
compactly supported weight $\varphi$. By shrinking $\sigma$ as necessary we may assume that $0 < \sigma \ll 1$.

Fix any compact time interval $[0,T]$, and let $K(T)$ be the best constant
for which the inequality
\begin{equation}\label{k-fuse}
 \int_0^T \| \langle x \rangle^{-1/2-\sigma} |\nabla e^{-itH} u_0|_g \|_{L^2(M)}^2
+ \| \langle x \rangle^{-3/2-\sigma} e^{-itH} u_0 \|_{L^2(M)}^2\ dt
\leq K(T)^2 \|u_0 \|_{\dot H^{1/2}(M)}^2
\end{equation}
holds for any Schwartz function $u_0$. From the local-in-time local smoothing theory in \cite{cks}, \cite{doi} we already 
know that $K(T)$ is finite for each $T$.
Our task is to show that $K(T)$ is bounded independently of $T$; the negative times can then be handled by time reversal symmetry. 

Recall that $M$ is equal to Euclidean space in the exterior region $|x| > R_0$.  To take advantage of this, let us fix $\varphi: M \to \R$
to be a smooth function which equals $1$ when $|x| \leq 4R_0$ and equals 0 when $|x| \geq 8R_0$.  We then define
the localized quantity $K_\varphi(T)$ to be the best constant such that
\begin{equation}\label{kvp-def}
 \int_0^T \| \varphi |\nabla e^{-itH} u_0|_g \|_{L^2(M)}^2\ dt
\leq K_\varphi(T)^2 \|u_0 \|_{\dot H^{1/2}(M)}^2
\end{equation}
holds for any Schwartz solution $u$ to \eqref{schrodinger}.  It is clear that $K_\varphi(T) \leq C_\varphi K(T)$.  In Section \ref{localize-sec} we shall establish the converse inequality
\begin{equation}\label{kvarphi}
K(T) \leq C_\varphi + C_\varphi K_\varphi(T)
\end{equation}
for all times $T > 0$, where the constants $C_\varphi$ depend on $\varphi$ but not on $T$.  

In light of \eqref{kvarphi}, we see that to bound $K(T)$ it suffices to bound the localized quantity $K_\varphi(T)$.  We make the
technical remark that $K_\varphi(T)$ is only required to control first derivatives of $e^{-itH} u_0$, and not $e^{-itH} u_0$ directly.
This will be important in the low freqency analysis later on.

The next step is energy decomposition into the very low energy, medium energy, and very high energy
portions of the evolution.  Let $0 < \eps_0 \ll 1$ be a small parameter (depending on $\varphi$) to be chosen later, and decompose
$1 = P_{lo} + P_{med} + P_{hi}$, where $P_{lo}$, $P_{med}$, $P_{hi}$ are the spectral multipliers
$$P_{lo} := \chi(H/\eps_0); \quad P_{med} = \chi(\eps_0 H) - \chi(H/\eps_0); \quad P_{hi} := 1 - \chi(\eps_0 H)$$
and $\chi: \R \to \R$ is a bump function supported on $[-1,1]$ which equals 1 on $[-1/2,1/2]$.  We shall prove the
following three propositions, in Sections \ref{lowfreq-sec}, \ref{medfreq-sec}, \ref{highfreq-sec} respectively:

\begin{proposition}[Low energy estimate]\label{lowfreq}  For any $T > 0$, we have
$$ \int_0^T \| \varphi \nabla P_{lo} e^{-itH} u_0 \|_{L^2(M)}^2\ dt \leq 
o_{\eps_0 \to 0; \varphi}(K(T)^2) \|u_0 \|_{\dot H^{1/2}(M)}^2.$$
\end{proposition}

\begin{proposition}[Medium energy estimate]\label{medfreq}  For any $T > 0$, any time-step $\tau \geq 1$ and radius $R \geq 10R_0$, 
we have
$$ \int_0^T \| \varphi \nabla P_{med} e^{-itH} u_0 \|_{L^2(M)}^2\ dt \leq 
(C_{\varphi,\eps_0,R,\tau} + o_{R \to \infty; \varphi, \eps_0}(K(T)^2) + o_{\tau \to \infty; \varphi,\eps_0,R}(K(T)^2)) \|u_0 \|_{\dot H^{1/2}(M)}^2.$$
\end{proposition}

\begin{proposition}[High energy estimate]\label{highfreq}  For any $T > 0$, we have
$$ \int_0^T \| \varphi \nabla P_{hi} e^{-itH} u_0 \|_{L^2(M)}^2\ dt \leq 
o_{\eps_0 \to 0; \varphi}(K(T)^2) \|u_0 \|_{\dot H^{1/2}(M)}^2.$$
\end{proposition}

Combining these three propositions using the triangle inequality and using the definition \eqref{k-fuse} of $K_\varphi(T)$, we see that
$$ K_\varphi(T)^2 \leq o_{\eps_0 \to 0,\varphi}(K(T)^2) + 
(C_{\varphi,\eps_0,R,\tau} + o_{R \to \infty; \varphi,\eps_0}(K(T)^2) + o_{\tau \to \infty; \varphi,\eps_0,R}(K(T)^2))
+ o_{\eps_0 \to 0;\varphi}(K(T)^2)$$
for any $\eps_0 > 0$, $\tau \geq 1$, and $R \geq 1$.
If we choose $\tau$ sufficiently large depending on $\varphi$, $\eps_0$ and $R$, and $R$ sufficiently large depending on $\varphi$, $\eps_0$, 
we conclude that
$$ K_\varphi(T)^2 \leq C_{\varphi,\eps_0} + o_{\eps_0 \to 0;\varphi}(K(T)^2);$$
combining this with \eqref{kvarphi} we obtain
$$ K(T)^2 \leq C_{\varphi,\eps_0} + o_{\eps_0 \to 0;\varphi}(K(T)^2);$$
letting $\eps_0$ be sufficiently small depending on $\varphi$ and recalling that $K(T)$ is finite, we conclude
that $K(T) \leq C_{\varphi}$ for all time $T$, which gives Theorem \ref{main}.

It remains to prove \eqref{kvarphi} and Propositions \ref{lowfreq}, \ref{medfreq}, \ref{highfreq}.  This will be done in the following sections.

\section{Physical space localization}\label{localize-sec}

We first prove \eqref{kvarphi}. 
Let $u = e^{-itH} u_0$ be a Schwartz solution to \eqref{schrodinger}.
Let $T > 0$.  We will allow all constants $C$ to depend on $\varphi$ and will no longer
mention this dependence explicitly. We normalize $\|u_0\|_{\dot H^{1/2}(M)} = 1$, which by unitarity of $e^{-itH}$ and
spectral calculus implies that
\begin{equation}\label{energy-est}
\sup_{t \in \R} \| u(t) \|_{\dot H^{1/2}(M)} = 1.
\end{equation}
Our task is to show that
$$
 \int_0^T \| \langle x \rangle^{-1/2-\sigma} \nabla u \|_{L^2(M)}^2
+ \| \langle x \rangle^{-3-\sigma} u \|_{L^2(M)}^2\ dt
\leq (C + C K_\varphi(T)^2) \| u_0 \|_{\dot H^{1/2}(M)}^2.$$
Note that it will not be relevant whether we measure the magnitude of $\nabla u$ using the metric $g$ or the Euclidean metric as they
only differ by at most a constant.

From \eqref{k-fuse} we already have
\begin{equation}\label{kf-1}
 \int_0^T \int_{|x| \leq 3R_0} |\nabla u|^2\ dx dt \leq C K_\varphi(T)^2.
\end{equation}
From an easy Poincar\'e inequality argument we also can show that
$$ \int_{|x| \leq 3R_0} |u|^2\ dx
\leq C (\int_{|x| \leq 4R_0} |\nabla u|^2\ dx + \int_{3R_0 \leq |x| \leq 4R_0} |u|^2\ dx.$$
Thus it will suffice to work in the Euclidean region $|x| > 3R_0$ and prove that
\begin{equation}\label{localized-targ}
\int_0^T \int_{|x| > 3R_0} |x|^{-3-2\sigma} |\nabla u|^2 + |x|^{-1-2\sigma} |u|^2\ dx dt \leq C + C K_\varphi(T)^2.
\end{equation}
We now invoke the positive commutator method.  Let $A$ be an arbitrary linear operator on Schwartz functions.  
From the self-adjoint nature of $H$, we observe the \emph{Heisenberg identity}
\begin{equation}\label{heisenberg}
 \frac{d}{dt} \langle Au(t), u(t) \rangle_{L^2(M)} = \langle i[H,A]u(t), u(t) \rangle_{L^2(M)}
\end{equation}
where $i[H,A] = i(HA-AH)$ is the Lie bracket of $H$ and $A$;
integrating this in $t$ and using \eqref{energy-est} and the duality of $\dot H^{1/2}(M)$ and $\dot H^{-1/2}(M)$, we obtain
\begin{equation}\label{commutator}
 |\int_0^T \langle i[H,A]u(t), u(t) \rangle_{L^2(M)}\ dx| \leq 
C \|A\|_{\dot H^{1/2}(M) \to \dot H^{-1/2}(M)}.
\end{equation}
The \emph{positive commutator method} is based on choosing $A$ so that
so that $i[H,A]$ is mostly positive definite, in order to extract useful information out of \eqref{commutator}.

Let us now set $A$ equal to the self-adjoint first-order operator
$$ A = -i a_{,k} \partial_k - i \partial_k a_{,k}$$
where $a: M \to \R$ is the function $a := \chi (|x| - \eps |x|^{1-\eps})$, where $\chi$ is a smooth cutoff supported on the region $|x| > 2R_0$ which equals 1 when $|x| \geq 3R_0$, and $0 < \eps \ll 1$ is a sufficiently small constant depending on $R_0$, $a_{,k}$ denotes the Euclidean derivative of $a$ in the $e_k$ direction, and we are summing indices in the usual manner.  Since $\nabla a = O(1)$ and $\nabla^2 a = O(1/|x|)$, we observe from \eqref{h1} and the classical Hardy inequality $\| u / |x| \|_{L^2(\R^3)} \leq C \| u \|_{\dot H^1(\R^3)}$ 
that $A$ maps $\dot H^1(M)$ to $L^2(M)$, and by self-adjointness also maps $L^2(M)$ to $\dot H^{-1}(M)$.  By interpolation
we conclude that $A$ maps $\dot H^{1/2}$ to $\dot H^{-1/2}$.
Also, since $H = -\frac{1}{2} \partial_j \partial_j$ on the support of $\chi$, we can compute
\begin{align*}
i[H,A] &= -\frac{1}{2} [\partial_j \partial_j, a_{,k} \partial_k + \partial_k a_{,k} ] \\
&= -\frac{1}{2} (a_{,kjj} \partial_k + 2 a_{,kj} \partial_{jk} + \partial_k a_{,kjj} + 2 \partial_k a_{,kj} \partial_j ) \\
&= - 2 \partial_j a_{,kj} \partial_k - \frac{1}{2} a_{,jjkk}
\end{align*}
and thus from \eqref{commutator} and an integration by parts we conclude that
\begin{equation}\label{x-ident}
|\int_0^T \int_{|x| > R_0} 
2 a_{,jk} \partial_k u \overline{\partial_j u} - \frac{1}{2} \Delta_{\R^3}^2 a |u|^2\ dx dt|
\leq C.
\end{equation}
Let us first consider the portion of \eqref{x-ident} on the region $|x| > 3R_0$, for which $a = |x| - \eps |x|^{1-\eps}$.  Then
a computation shows that $a$ has some convexity if $\eps$ is sufficiently
small, indeed in this region we have
$$ a_{,jk} \partial_k u \overline{\partial_j u} \geq c_\eps \frac{|\nabla u|^2}{|x|^{1+\eps}}; \quad \Delta^2 a \leq - c_\eps / |x|^{3+\eps}$$
for some $c_\eps > 0$.
Invoking \eqref{x-ident} and using \eqref{kf-1} to estimate the region where $2R_0 \leq |x| \leq 3R_0$, we conclude that
$$
2 \int_0^T \int_{|x| > 3R_0} \frac{|\nabla u|^2}{|x|^{1+\eps}} + c_\eps \frac{|u|^2}{|x|^{2+\eps}}\ dx dt 
\leq C_\eps + C_\eps K_\varphi(T)^2 + \frac{1}{2} \int_0^T \int_{|x| \leq 3R_0} \Delta_{\R^3}^2 a |u|^2\ dx dt.
$$
To conclude \eqref{localized-targ}, it thus suffices by \eqref{kf-1} to establish following the fixed-time estimate:

\begin{lemma}\label{fixed} If $\eps > 0$ is sufficiently small, then there exists a constant $C > 0$ such that
$$ \int_{|x| \leq 3R_0} \Delta_{\R^3}^2 a |f|^2 \leq C \int_{|x| \leq 3R_0} |\nabla f|^2$$
for all smooth functions $f$.  (Note that the left-hand side can be negative).
\end{lemma}

\begin{proof} It is possible to establish this from Poincar\'e inequality and the Green's function computation below (the main point being
that $\Delta^2 a$ has negative mean), but we shall use a compactness argument instead.
Let $\delta > 0$ be a small number to be chosen later.  It clearly suffices to show that there exists $C > 0$
such that
$$ \int_{|x| \leq 3R_0} (\delta + \Delta_{\R^3}^2 a) |f|^2 \leq \delta \int_{|x| \leq 3R_0} |f|^2+
C \int_{|x| \leq 3R_0} |\nabla f|^2$$
for all smooth functions $f$.  We have chosen to use the Euclidean measure $dx$ here but one could equally well run the following argument
using the measure $dg$.

Suppose for contradiction that the above estimate failed.  Then we can find a sequence $f_n$ of smooth functions with the normalization
$$ \delta \int_{|x| \leq 3R_0} |f_n|^2 + n \int_{|x| \leq 3R_0} |\nabla f_n|^2 = 1$$
such that
$$ \limsup_{n \to \infty} \int_{|x| \leq 3R_0} (\delta + \Delta^2 a) |f_n|^2 > 1.$$
By Rellich compactness we can find a subsequence $f_{n_j}$ of $f_n$ which converges in $L^2$ to a limiting object $f \in L^2(B(0,3R_0))$, and 
then by Fatou's lemma
$$ \int_{|x| \leq 3R_0} |\nabla f|^2 = 0.$$
In other words, $f$ is equal to a constant on the region $\langle x \rangle \leq 3R_0$.  But by Green's formula we have
$$ \int_{|x| \leq 3R_0} \delta + \Delta^2 a = O(\delta) + \int_{|x| = 3R_0} \frac{d}{dn} \Delta a\ dS = 
O(\delta) + 4\pi (3R_0)^2 (-\frac{1}{(3R_0)^2} + O(\eps)) = 
-4\pi + O(\eps) + O(\delta)$$
which is negative if $\eps$ and $\delta$ are chosen sufficiently small.  This is a contradiction, and the claim follows.
\end{proof}

\section{Low energy estimate}\label{lowfreq-sec}

Now we prove the low energy estimate in Proposition \ref{lowfreq}, which is the easiest of the three propositions to prove, especially
in the model case when there is no potential $V$, and the manifold $M$ is a compact perturbation of $\R^3$.  The idea here is
to exploit the uncertainty principle to extract some gain from the spatial projection $\varphi$ and the frequency
projection $P_{lo}$, but in order to do this we need to somehow exploit the fact that $H$ contains no resonances or bound states at the zero energy.
For technical reasons (having to do with our use of the homogeneous norm $\dot H^{1/2}(M)$ instead of the inhomogeneous norm $H^{1/2}(M)$)
we shall need the following non-standard formulation of this non-resonance property.

\begin{proposition}[Laplace-Beltrami operators have no resonance]\label{nonresonance}  Suppose that $f: M \to \C$ is a measurable, weakly differentiable function such that
$$
 \int_M \langle x \rangle^{-3-2\sigma} |f(x)|^2 + \langle x \rangle^{-1-2\sigma} |\nabla f|^2 \ dg(x) < \infty
$$
for some $K < \infty$, and such that $Hf \equiv 0$ in the sense of distributions.  Then (if $\sigma > 0$ is sufficiently
small) $f$ is a  constant.  
\end{proposition}

\begin{proof} We may take $f$ to be real.  The condition $Hf=0$, combined with the local square-integrability of $f$, implies that $f$ is in
fact smooth thanks to elliptic regularity (and the smoothness of $g$).  
Since $Hf=0$ and $H = H_0 = -\frac{1}{2} \Delta_{\R^3}$ outside of $B(0,R_0)$, the function 
$\Delta_{\R^3} f$ is a smooth compactly supported function.   Thus if we set $F := \frac{1}{4\pi |x|} * \Delta_{\R^3} f$, i.e. the convolution of $\Delta_{\R^3} f$ with the fundamental solution of the Euclidean Laplacian, then $f-F$ is harmonic, $F$ decays like $O( 1/ \langle x \rangle )$, and $|\nabla F|$ decays like $O(1/\langle x \rangle^2)$.  From hypothesis and the triangle inequality we then have
$$ \int_M \langle x \rangle^{-3-2\sigma} |f(x)- F(x)|^2\ dg(x) < \infty.$$
Since $f-F$ is harmonic, we conclude (e.g. from the mean-value theorem applied to $f-F$ and its first derivative) that $f-F$ is constant,
thus $f-F = c$.  By subtracting this constant from $f$ we may in fact take $F=f$.  This now shows that $f(x) = O(1/\langle x \rangle)$
and $\nabla f(x) = O(1/\langle x\rangle^2)$.  But then we can justify the computation
$$ \int_M |\nabla f|_g^2\ dx = - \int_M f(x) \nabla^j \nabla_j f(x)\ dg(x) = \int_M 2 f(x) Hf(x)\ dg(x) = 0$$
by inserting a suitable smooth cutoff to a large ball $B(0,R)$ and then letting $R \to \infty$; we omit the standard details.  But
then we have $\nabla f = 0$ and we are done.
\end{proof}

Using this fact and a compactness argument, we now conclude

\begin{proposition}[Poincar\'e-type inequality]\label{poincare} Let $f: M \to \R$ be a Schwartz function.  Let $\varphi$
be as in previous sections.  Then for any $\eps > 0$ we have
\begin{equation}\label{poincare-3}
\begin{split}
 \| \varphi |\nabla f| \|_{L^2(M)}
\leq o_{\eps \to 0}(1)
(&\| \langle  x \rangle^{-1/2-\sigma} \nabla f \|_{L^2(M)} 
+
\| \langle  x \rangle^{-3/2-\sigma} f \|_{L^2(M)} \\
&+ \eps^{-1} \| \langle  x \rangle^{-1/2-\sigma} \nabla H f \|_{L^2(M)}
+ \eps^{-1} \| \langle  x \rangle^{-3/2-\sigma} H f \|_{L^2(M)}).
\end{split}
\end{equation}
\end{proposition}

\begin{proof}  Suppose for contradiction that the Proposition was false.  Then there exists a $\delta > 0$,
a sequence $\eps_n > 0$ converging to zero, and Schwartz functions $f_n$ such that
\begin{align*}
 \| \varphi |\nabla f_n| \|_{L^2(M)} > \delta
(&\| \langle  x \rangle^{-1/2-\sigma} \nabla f_n \|_{L^2(M)} +
\| \langle  x \rangle^{-3/2-\sigma} f_n \|_{L^2(M)}\\
&+ \eps_n^{-1} \| \langle  x \rangle^{-1/2-\sigma} \nabla H f_n \|_{L^2(M)}
+ \eps_n^{-1} \| \langle  x \rangle^{-3/2-\sigma} H f_n \|_{L^2(M)}).
\end{align*}
Without loss of generality we may assume that $f_n$ is not identically zero, and then we can normalize
the expression in parentheses to equal 1.  Thus
\begin{equation}\label{compactness}
\begin{split}
 \| \langle  x \rangle^{-1/2-\sigma} \nabla f_n \|_{L^2(M)},
\| \langle  x \rangle^{-3/2-\sigma} f_n \|_{L^2(M)} &\leq 1; \\
\| \langle  x \rangle^{-1/2-\sigma} \nabla H f_n \|_{L^2(M)} + 
\| \langle  x \rangle^{-3/2-\sigma} H f_n \|_{L^2(M)} &\leq \eps_n; \\
\quad  \| \varphi |\nabla f| \|_{L^2(M)} &\geq \delta.
\end{split}
\end{equation}
Next, we establish weighted
$\dot H^2$ bounds on $f_n$
via a Bochner identity.  From \eqref{compactness} and Cauchy-Schwarz
we have
$$ \int_M \langle x \rangle^{-1-2\sigma} \Re((\nabla^j f_n) \overline{\nabla_j H f_n}) \leq \eps_n = O(1).$$
We substitute $H = - \frac{1}{2} \nabla^k \nabla_k$.  Since $M$ is flat outside of the compact set $B(0,R_0)$ we have
$$ \nabla_j H f_n = -\frac{1}{2} \nabla_j \nabla^k \nabla_k f_n = -\frac{1}{2} \nabla^k \nabla_j \nabla_k f_n + 
O( |\nabla f_n|) \chi_{B(0,R_0)}.$$
Using this and \eqref{compactness}, we obtain
$$ -\int_M \langle x \rangle^{1-2\sigma} \Re( (\nabla^j f_n) \nabla^k \nabla_j \nabla_k \overline{f_n} ) \leq C.$$
Integrating by parts we obtain
$$ \int_M \langle x \rangle^{1-2\sigma} \Re( (\nabla^k \nabla^j f_n) \nabla_j \nabla_k \overline{f_n} ) 
+ (\nabla^k \langle x \rangle^{1-2\sigma}) \Re( (\nabla^j f_n) \nabla^k \nabla_j \nabla_k \overline{f_n} ) \leq C$$
Since $\Re( (\nabla^j f_n) \nabla^k \nabla_j \nabla_k \overline{f_n} ) = \frac{1}{2} \nabla_k
|\nabla f_n|_g^2$, we can integrate by parts once more to obtain
$$ \int_M \langle x \rangle^{1-2\sigma} |\Hess(f_n)|_g^2 \leq C+ \frac{1}{2} \int_M 
(\Delta_M \langle x \rangle^{1-2\sigma}) |\nabla f_n|_g^2.$$
Since $\Delta_M \langle x \rangle^{1-2\sigma} = O(\langle x \rangle^{-1-2\sigma})$, 
the integral on the right-hand side is $O(1)$ by \eqref{compactness}.  
Thus
\begin{equation}\label{hessian}
  \int_M \langle x \rangle^{1-2\sigma} |\Hess(f_n)|^2 \leq C.
\end{equation}
From this, \eqref{compactness} and Rellich compactness we see that the sequence $f_n$, when localized smoothly to
any ball $B(0,R)$, is contained in a compact subset of $H^1(B(0,R))$.  From this and the usual Arzela-Ascoli diagonalization
argument, we may extract a subsequence of $f_n$ which converges locally in $H^1$ to
some limit $f$.  From \eqref{compactness} and Fatou's lemma we then see that
\begin{equation}\label{compactness-final}
\| \langle  x \rangle^{-1/2-\sigma} \nabla f \|_{L^2(M)},
\| \langle  x \rangle^{-3/2-\sigma} f \|_{L^2(M)} \leq 1; \quad \| \varphi |\nabla f| \|_{L^2(M)} \geq \delta
\end{equation}
and that $Hf=0$ in the sense of distributions.   But then by Proposition \ref{nonresonance}, $f$ is constant.
But this contradicts the last estimate in \eqref{compactness-final}, and we are done.
\end{proof}

We can now quickly prove Proposition \ref{lowfreq}. Applying \eqref{k-fuse} with $u_0$ replaced by $P_{lo} u_0$ and
$H P_{lo} u_0$, and using some spectral theory to estimate the right-hand side, we obtain
$$ \int_0^T \| \langle x \rangle^{-1/2-\sigma} \nabla P_{lo} e^{-itH} u_0 \|_{L^2(M)}^2
+ \| \langle x \rangle^{-3/2-\sigma} P_{lo} e^{-itH} u_0 \|_{L^2(M)}^2\ dt
\leq C K(T)^2 \|u_0 \|_{\dot H^{1/2}(M)}^2$$
and
$$ \int_0^T \| \langle x \rangle^{-1/2-\sigma} \nabla H P_{lo} e^{-itH} u_0 \|_{L^2(M)}^2
+ \| \langle x \rangle^{-3/2-\sigma} P_{lo} e^{-itH} u_0 \|_{L^2(M)}^2\ dt
\leq C \eps_0^{2} K(T)^2 \|u_0 \|_{\dot H^{1/2}(M)}^2;$$
note that $P_{lo}$ and $H$ are spectral multipliers and hence commute with each other and with $e^{-itH}$.
Applying Proposition \ref{poincare}, Proposition \ref{lowfreq} follows. 

\section{High energy estimate}\label{highfreq-sec}

We now prove Proposition \ref{highfreq}, which is the next easiest of the three propositions.  This case resembles
the local-in-time theory of Craig-Kappeler-Strauss \cite{cks} and Doi \cite{doi}, and indeed our main tool here will 
be the positive commutator method applied to a certain pseudo-differential operator, exploiting the non-trapping 
hypothesis to ensure that the symbol of the pseudo-differential operator increases along geodesic flow.
As we shall now be working in the high energy setting, we will not need to take as much care with lower order terms
as in previous sections.  For similar reasons, we will not need to use the homogeneous Sobolev spaces $\dot H^s(M)$, relying instead
on the more standard (and more stable) inhomogeneous Sobolev spaces $H^s(M)$.  The argument here is in fact quite general and would work on
any asymptotically conic manifold with a short-range metric perturbation and a short-range potential.

It will be convenient
to use the \emph{scattering pseudo-differential calculus}, which is an extension of the standard pseudo-differential calculus which keeps track of the decay of the symbol at infinity.  We briefly summarize
the relevant features of this calculus here, referring the reader to \cite{cks} for more complete details. 
For any $m,l \in \R$, we define a \emph{symbol $a: T^* M \to \C$ of order $(m,l)$} to be any
smooth function obeying the bounds
$$ |\nabla_x^\alpha \nabla_\xi^\beta a(x,\xi)| \leq 
C_{\alpha,\beta} \langle \xi \rangle^{m-|\beta|} \langle x \rangle^{-l-|\alpha|};$$
the function $a(x,\xi) = \langle x \rangle^{-l} \langle \xi \rangle^m$ is a typical example of such a symbol.
Note that we assume that each derivative in $x$ gains a power of $\langle x \rangle$,
in contrast to the standard symbol calculus in which no such gain is assumed.  
We let $S^{m,l}(\Mbar)$ denote the space of such symbols.  Given any such symbol $a \in S^{m,l}(\Mbar)$, 
we can define an associated pseudo-differential operator $A = \Op(a)$ by the usual 
Kohn-Nirenberg quantization formula
$$
\Op(a) u(x) := (2\pi)^{-n} \int e^{i\langle x-y, \xi \rangle} a(x,\xi) u(y)\, dy \, d\xi.
$$
We sometimes denote $a$ by $\sigma(A)$ and refer to it as the \emph{symbol} of $A$.  Heuristically speaking,
we have $A = \sigma(A)(x, \frac{1}{i} \nabla_x)$.  We refer to the class of pseudo-differential operators
of order $(m,l)$ as $\Psisc^{m,l}$.
Also, if $h: \R \to \C$ is any spectral symbol of order $m/2$, the corresponding spectral multiplier $h(H)$ is
a pseudo-differential operator of order $(m, 0)$.  In particular, 
$(1+H)^{m/2}$ has order $(m,0)$, and the Littlewood-Paley type 
operators $P_{lo}$, $P_{med}$, $P_{hi}$ have order $(0,0)$.  We caution however that the Schr\"odinger 
propagators $e^{-itH}$ are not pseudo-differential operators.

The composition of an operator $A = \Op(a)$ of order $(m,l)$ with an operator of $B = \Op(b)$ order $(m',l')$ is an
operator $AB$ of order $(m+m', l+l')$, whose symbol $\sigma(AB)$ is equal to $\sigma(A) \sigma(B)$ plus an error
of order $(m+m'-1,l+l'+1)$; note the additional gain of 1 in the decay index $l$, which is not present in the
classical calculus.  
Similarly, the commutator $i[A,B]$ will be an operator of order $(m+m'-1,l+l'-1)$ with
symbol $\sigma(i[A,B])$ equal to the Poisson bracket
$$ \{ \sigma(A), \sigma(B) \} := \nabla_x \sigma(A) \cdot \nabla_\xi \sigma(B) - \nabla_\xi \sigma(A) \cdot \nabla_x \sigma(A),$$
plus an error of order $(m+m'-2, l+l'+2)$. We shall write the above facts schematically as
$$ \sigma(AB) = \sigma(A) \sigma(B) + O( S^{m+m'-1,l+l'+1});
\quad \sigma( i[A,B] ) = \{ \sigma(A), \sigma(B) \} + O( S^{m+m'-2,l+l'+2} )$$
or equivalently as
$$ \Op(a) \Op(b) = \Op(ab) + O( \Psisc^{m+m'-1,l+l'+1}),
\quad i[\Op(a), \Op(b)] = \Op( \{a,b\} ) + O( \Psisc^{m+m'-2,l+l'+2}).$$
In particular, since $H$ has order $(2,0)$ and has principal symbol $\frac{1}{2} |\xi|_{g(x)}^2$
plus lower order terms of order $(1,1)$ and $(0,2)$, we see that if $a \in S^{m,l}$, then we have
$$ i[H, \Op(a)] = \Op( X a ) + O( \Psisc^{m,l+2}),$$
where $X a$ denotes the derivative of $a$ along geodesic flow in the cotagent bundle $T^* M$.  

Associated with the scattering calculus are the weighted Sobolev spaces $H^{m,l}(M)$ defined (for instance)
by
$$ \| u \|_{H^{m,l}(M)} := \| \langle x \rangle^l (1 + H)^{m/2} u \|_{L^2(M)}$$
(many other equivalent expressions for this norm exist, of course); when $l=0$ this corresponds to the usual Sobolev space $H^m(M)$.  It is easy to verify that a
scattering pseudo-differential operator of order $(m,l)$ maps $H^{m',l'}(M)$ to $H^{m'-m,l'+l}(M)$ for any $m', l'$.

In \cite{cks} (see also \cite{doi}) it was shown that the non-trapping hypothesis on $M$ allows one to
construct a real-valued symbol $a \in S^{1,0}$ (depending on $\varphi$) which was non-decreasing along
geodesic flow, $Xa \geq 0$, and in fact obeyed the more quantitative estimate
$$ Xa(x,\xi) = \varphi(x) |\xi|_g^2 + |b|^2$$
for some symbol $b$ of order $(1,1/2-\sigma)$.  In Euclidean space, an example of such an symbol is
$C_\varphi \frac{x}{\langle x \rangle} \cdot \xi$ for some sufficiently large constant $C_\varphi$.  Quantizing
this, we obtain
$$ i[H,A] = \nabla^j \varphi(x) \nabla_j + B^* B + O( \Psisc^{1,2-2\sigma} ),$$
where $A := \Op(a)$ is a symbol of order $(1,0)$, and $B := \Op(b)$ is a symbol of order $(1,1/2-\sigma)$.
We then apply the self-adjoint projection $P_{hi} = P^*_{hi}$ to both sides, 
and observe that this commutes with $H$, to obtain
$$ i[H,P^*_{hi} A P_{hi}] = P^*_{hi} \nabla^j \varphi(x) \nabla_j P_{hi} + 
P_{hi}^* B^* B P_{hi} + P_{hi} O( \Psisc^{1,2-2\sigma} ) P_{hi}.$$
Applying \eqref{commutator}, and integrating by parts (discarding the positive term $B^* B$, and using
that $\Psisc^{1,2-2\sigma}$ maps $H^{1/2,-1+\sigma}$ to $H^{-1/2,1-\sigma}$) we obtain
$$
 \int_0^T \int_M \varphi |\nabla e^{-itH} u_0|_g^2\ dg dt \leq 
C \|P^*_{hi} A P_{hi}\|_{\dot H^{1/2}(M) \to \dot H^{-1/2}(M)} \| u_0 \|_{\dot H^{1/2}(M)}^2
+ C \int_0^T \| P_{hi} e^{-itH} u_0 \|_{H^{1/2,-1+\sigma}}^2 \ dt.
$$
Since $A$ is of order $(1,0)$, and $P_{hi}$ maps the homogeneous Sobolev spaces to their inhomogeneous counterparts,
we see that $\|P^*_{hi} A P_{hi}\|_{\dot H^{1/2}(M) \to \dot H^{-1/2}(M)}$ is bounded by some $C_\varphi$.  
To finish the proof of
Proposition \ref{highfreq}, it thus suffices to show that
$$ \int_0^T \| P_{hi} e^{-itH} u_0 \|_{H^{1/2,-1+\sigma}}^2 \ dt
\leq o_{\eps_0}(1) \|u(0) \|_{\dot H^{1/2}(M)}^2.$$
On the other hand, applying \eqref{k-fuse} to $H^{-j} P_{hi} u(0)$ for $j=0,1$
 we have
$$
 \int_0^T \| \nabla H^{-j} e^{-itH} P_{hi} u_0 \|_{H^{0,-1/2-\sigma}}^2 
+ \| H^{-j} e^{-itH} P_{hi} u_0 \|_{H^{0,-3/2-\sigma}}^2 
\ dt
\leq C K(T) ^2\|u(0) \|^2_{H^{1/2}(M)}$$
for $j=0,1$. Thus fixing $t$ and setting $f := H^{-1} e^{-itH} P_{hi} u_0$
it will suffice to prove the fixed-time estimate
$$ \| H f \|_{H^{1/2,-1+\sigma}(M)}
\leq o_{\eps_0\to 0}(1) \sum_{j=0}^1\Big ( \| \nabla H^j f \|_{H^{0,-1/2-\sigma}(M)}
+ \| H^j f \|_{H^{0,-3/2-\sigma}(M)}\Big ).$$
It suffices to verify this for real-valued $f$.
By Rellich compactness, the space $H^{1, -1/2-\sigma}(M)$ embeds compactly into
$H^{1/2,-1+\sigma}(M)$.  Since $P_{hi} f \to 0$ as $\eps_0 \to 0$ for each individual $f$, 
we thus see by compactness that it suffices to show that
\begin{align}
 \| H f \|_{H^{1,-1/2-\sigma}(M)}
&\leq C ( \| \nabla H f \|_{H^{0,-1/2-\sigma}(M)} +
\| H f \|_{H^{0,-3/2-\sigma}(M)} \nonumber \\&+
\| \nabla f \|_{H^{0,-1/2-\sigma}(M)}
+ \| f \|_{H^{0,-3/2-\sigma}(M)} ).\label{eq:a-H}
\end{align}
The top order term of $\| H f \|_{H^{1,-1/2-\sigma}(M)}$ is already controlled by the right-hand side of \eqref{eq:a-H}, so
it suffices to control the lower order term $\| Hf \|_{H^{0,-1/2-\sigma}(M)}$.  But an integration by
parts allows us to write
$$ \int_M \langle x \rangle^{-1-2\sigma} Hf Hf\ dg = 
\frac{1}{2} \int_M \langle x \rangle^{-1-2\sigma} \Big (\nabla^j f \nabla_j Hf\ dg - 
(1+2\sigma) \frac {\nabla^j \langle x\rangle}{\langle x\rangle} \nabla^j f  Hf\Big )
.$$
The claim then follows from the Cauchy-Schwarz inequality.
This concludes the proof of Proposition \ref{highfreq}.

\section{Medium energy estimate}\label{medfreq-sec}

We now turn to the medium energy estimate, Proposition \ref{medfreq}, which is the hardest of the three propositions.
Neither the uncertainty principle nor the non-trapping condition will be much use here.  Instead our tools\footnote{One can also proceed via Kato's theory of $H$-smooth operators, using the limiting absorption principles obtained in \cite{burq}.  We will pursue this approach in detail in \cite{rt:limiting}.} will be
a RAGE-type theorem, exploiting the fact that $H$ has no embedded eigenvalues, to propagate the solution away from
the origin after a long time $\tau$ (though the energy localization shows the solution will not move \emph{too} far away from the origin in bounded time), combined with a decomposition of phase space into incoming and 
outgoing waves (cf. Enss' method, \cite{enss}, \cite{simon1}), and several applications of Duhamel's formula.  Because we have eliminated the high frequencies,
we will enjoy an approximate finite speed of propagation law for the solution (but the upper bound of the speed is quite
large, being roughly $O(\eps_0^{-1/2})$), and because we have eliminated the low energies, we will not encounter
a frequency singularity when we decompose into incoming and outgoing waves, although again we will pick up some negative
powers of $\eps_0$.  We shall compensate for these $\eps_0$ losses by using the RAGE theorem to gain a factor of
$o_{\tau \to \infty; \eps_0, R}(1)$ within a distance $R$ from the origin, and to also gain a factor of
$o_{R \to \infty; \eps_0}(1)$ when one is further than $R$ from the origin.  The reader should view the comparative
magnitudes of $\eps_0, R, \tau$ according to the relationship $1/\eps_0 \ll R \ll \tau$, which is of course
the most interesting case of Proposition \ref{medfreq}.

It is instructive at this point to  recall
the basic features of the Enss' method. Define the wave operator
$W_+=s-\lim_{t\to +\infty}  e^{itH} e^{-itH_0}$. Enss' method is concerned with establishing
{\it completeness}
of $W_{+}$, i.e., showing that the range of $W_{+}$ coincides with the 
continuous subspace of $H=H_0+V$. We assume otherwise so that there exists 
$\phi_0$ not in the range of $W_+$. Density arguments allow us to consider 
$\phi_0$ with compact support and medium energies only. We then evolve $\phi_0$
by $e^{-itH}$ and claim that we can find a sequence of times $t_n$ and a decomposition
\begin{equation}\label{ndecomp}
e^{-it_nH}\phi_0=\phi_n=\phi_{n,out} +\phi_{n,in} + \tilde \phi_{n},
\end{equation}
(compare this with the decomposition into $F_{loc}, F_{glob}$ in \eqref{flg-def} and 
a further decomposition of $F_{glob}$ into the outgoing/incoming waves in 
\eqref{in-and-out})
with the properties that 
$$\|(W_{+}-1)\phi_{n,out}\|_{L^2}, \|\phi_{n,in}\|_{L^2}, \|\tilde \phi_n\|_{L^2} = o_{n \to \infty}(1).$$
If we had such a decomposition, we could conclude from the $L^2$ boundedness of $W_+$ that
$$
\|W_+ e^{it_nH_0} \phi_0-\phi_0\|_{L^2} = o_{n \to \infty}(1)
$$
and thus $\phi_0$ lies in the range of $W_+$. The desired
decomposition \eqref{ndecomp} can be found for instance in \cite{simon1}, and we sketch it as follows. 
The local component $\tilde \phi_n$ can be set for instance to $\tilde \phi_n := \chi_{|x| \leq n} \phi_n$,
its convergence to zero is a consequence of the RAGE theorem (if the times $t_n$ are chosen appropriately). 
The global component $(1-\chi_{|x|\le n} ) \phi_n$ is partitioned into functions $\phi_{n,\alpha}$ 
supported in unit balls with centers at the lattice points $\alpha\in \R^n$. One can then define
the 
$$
\phi_{n,\alpha,out}=\int_{\xi,y} e^{i(x-y)\cdot \xi} m(\xi) \phi_{n,\alpha}(y)\,d\xi\,dy,
\qquad \phi_{n,\alpha,in}=
\int_{\xi,y} e^{i(x-y)\cdot \xi} (1-m(\xi)) \phi_{n,\alpha}(y)\,d\xi\,dy,
$$   
where $m(\xi)$ is a smooth multiplier localizing to the region
$\angle (\xi,\alpha)\le 3\pi/2$. Since it is expected that a compact support 
function $\phi_0$
propagated by $e^{-it_n H}$ should become mostly outgoing in the region 
$|x|\ge n$ we have that $\phi_{n,in}\to 0$,  while the fact that 
$(W_+-1)\phi_{n,out}\to 0$ in $L^2$ follows since on the outgoing waves 
the evolution $e^{-itH}$ is well approximated by the free flow $e^{-itH_0}$ 
and in fact converges to it as $n\to \infty$.

Our proof of Proposition \ref{medfreq} will follow in spirit the above construction, although we have 
to take additional care since we do not work with functions of compact support and  
need weighted estimates instead of $L^2$ bounds.  On the other hand, we are still localized in energy,
and so we will not be too concerned about losing or gaining too many derivatives (as we are able to lose factors
of $\eps_0$ or $\eps_0^{-1}$ in our estimates here).  
In particular, the metric perturbation $H$ of $H_0$ is now
of similar ``strength'' to a potential perturbation, and we will now be able to use $H_0$ as a reasonable approximant
to $H$, at least when the solution is far away from the origin.

It is more convenient to work in the dual formulation.  First observe that the claim is easy when $T \leq \tau$,
since $P_{med} e^{-itH}$ maps $\dot H^{1/2}(M)$ to $H^1(M)$ thanks to the frequency localization of $P_{med}$.
Similarly the $\int_0^{\tau}$ portion of the integral is easy to deal with.  Thus we may assume that $T > \tau$
and reduce to proving
$$ \int_{\tau}^T \| \varphi \nabla P_{med} e^{-itH} u_0 \|_{L^2(M)}^2\ dt \leq 
(C_{\eps_0,R,\tau} + o_{R \to \infty; \eps_0}(K(T)) + o_{\tau \to \infty; \eps_0,R}(K(T))) \|u_0 \|_{\dot H^{1/2}(M)}^2,$$
which after dualization (and shifting $t$ by $\tau$) becomes
\begin{equation}\label{dual}
\begin{split}
 \| \int_0^{T-\tau} e^{itH} e^{i\tau H} P_{med} \nabla_j \varphi F^j(t+\tau)\ dt \|_{\dot H^{-1/2}(M)}
\leq &
(C_{\eps_0,R,\tau} + o_{R \to \infty; \eps_0}(K(T)) + o_{\tau \to \infty; \eps_0,R}(K(T)))\\ 
&(\int_0^{T-\tau} \| F(t+\tau) \|_{L^2(M)}^2\ dt)^{1/2}
\end{split}
\end{equation}
for any vector field $F(t)$ defined on the time interval $[\tau,T]$ which is Schwartz for each time $t$.

Fix $F$. We split
\begin{equation}\label{decomp}
 e^{itH} e^{i\tau H} P_{med} \nabla_j \varphi F^j(t+\tau) = e^{itH} F_{loc}(t) + e^{itH} F_{glob}(t)
\end{equation}
where
\begin{equation}\label{flg-def}
F_{loc}(t) := \varphi_R e^{i\tau H} P_{med} \nabla_j \varphi F^j(t+\tau); \quad
F_{glob}(t) := (1-\varphi_R) e^{i\tau H} P_{med} \nabla_j \varphi F^j(t+\tau).
\end{equation}

\divider{Term 1: Contribution of the local part.}

Consider first the contribution of the local term $F_{loc}(t)$. To control this term we use the following local decay result.

\begin{proposition}[RAGE theorem]\label{rage} Let $\varphi_R$ be a bump function supported on $B(0,2R)$ which equals 1 on $B(0,R)$. For all Schwartz vector fields $f^j$, we have
$$ \| \varphi_R e^{i\tau H} P_{med} \nabla_j \varphi f^j \|_{L^2(M)}
\leq o_{\tau \to \infty; R, \eps_0}( \| f \|_{L^2(M)} ).$$
\end{proposition}

\begin{proof} Observe that we have the crude bound
$$ \| \varphi_R e^{i\tau H} P_{med} \nabla_j \varphi f^j \|_{L^2(M)} \leq C_{\eps_0} \| f \|_{H^{-1,-1}(M)}$$
since $P_{med} \nabla \varphi$ maps $H^{-1,-1}$ to $L^2$, and $\varphi_R e^{i\tau H}$ is bounded on $L^2$.
Since the unit ball of $L^2(M)$ is precompact in $H^{-1,-1}$, it thus suffices 
to prove the estimate
$$ \| \varphi_R e^{i\tau H} P_{med} \nabla_j \varphi f^j \|_{L^2(M)}
\leq o_{\tau \to \infty; R, f, \eps_0}( 1 )$$
for all Schwartz vector fields $f$; the point being that the compactness allows us to ignore the dependence of the constants
on $f$.  
Fix $\tilde f$, $\varphi$, $R$.  Since $f$ is Schwartz, the 
curve $\{ \tilde \varphi e^{i\tau H} P_{med} \nabla_j \varphi f^j: \tau \in \R \}$ is bounded in $H^{1,1}(M)$ (for instance) and hence precompact in $L^2$.  Thus to prove the above strong 
convergence, it will suffice to prove the weak convergence result
$$ |\langle \tilde \varphi e^{i\tau H} P_{med} \nabla_j \varphi f^j, \psi \rangle|
\leq o_{\tau \to \infty; R, f, \eps_0, \psi}( 1 )$$
for all Schwartz functions $\psi$.

Fix $\psi$.  Since the spectrum of $H$ is purely absolutely continuous\footnote{Actually, the argument would still work well if $H$ had some singular continuous spectrum, except that $\tau$ must now be averaged over an interval such as $[\tau_0, 2\tau_0]$, but the reader may verify that the arguments below will continue to work with this averaging.  The spectral fact which is really being used here is that $H$ contains no embedded eigenfunctions at medium frequencies, since such eigenfunctions would certainly contradict \eqref{ls-manifold-global}.}, we see from the Riemann-Lebesgue lemma and
the spectral theorem that
$$ \langle \varphi_R e^{i\tau H} P_{med} \nabla \varphi f, \psi \rangle = \langle e^{i\tau H} P_{med} \nabla_j \varphi f^j, \varphi_R \psi \rangle
= o_{\tau \to \infty; P_{med} \nabla_j \varphi f^j, \varphi_R \psi}(1)$$
and the claim follows.
\end{proof}

From this proposition we see in particular that
$$ \| F_{loc}(t) \|_{L^2(M)} \leq o_{\tau \to \infty; R, \eps_0}( 1 ) \| F(t+\tau) \|_{L^2(M)}.$$
From dualizing the second part of \eqref{k-fuse} we have
\begin{equation}\label{Kdual}
 \| \int_0^T e^{itH} \langle x \rangle^{-3/2-\sigma} G(t)\ dt \|_{\dot H^{-1/2}(M)}
\leq K(T) (\int_0^T \| G(t) \|_{L^2(M)}^2\ dt)^{1/2}
\end{equation}
for any $G$.
If we truncate the time interval to $T-\tau$, substitute $G(t) := \langle x \rangle^{3/2+\sigma} F_{loc}(t)$,
and take advantage of the spatial localization of $F_{loc}$ to $B(0,2R)$, we thus have
$$ \| \int_0^{T-\tau} e^{itH} F_{loc}(t)\ dt \|_{\dot H^{-1/2}(M)}
\leq o_{\tau \to \infty; \eps_0, R}( K(T) ) (\int_0^{T-\tau} \| F(t+\tau)\|_{L^2(M)}^2\ dt)^{1/2}.$$
Thus the first term in \eqref{decomp} is acceptable.  

\divider{Term 2: Contribution of the global part.}

To conclude the proof of Proposition \ref{medfreq}, it suffices
to establish the estimate
\begin{equation}\label{glob-duhamel}
 \| \int_0^{T-\tau} e^{itH} F_{glob}(t)\ dt \|_{\dot H^{-1/2}(M)}
\leq (C_{\eps_0,R,\tau} + o_{R \to \infty; \eps_0}( K(T) )) (\int_0^{T-\tau} \| F(t+\tau)\|_{L^2(M)}^2\ dt)^{1/2}
\end{equation}
for the global term $F_{glob}(t)$.

Let us first show an associated estimate for the \emph{free} propagator $e^{itH_0}$, namely
\begin{equation}\label{glob-free}
 \| \int_0^{T-\tau} e^{itH_0} F_{glob}(t)\ dt \|_{\dot H^{-1/2}(M)}
\leq C_{\eps_0,R,\tau} (\int_0^{T-\tau} \| F(t+\tau)\|_{L^2(M)}^2\ dt)^{1/2}.
\end{equation}
Note that our constants are allowed to depend on $\tau$, as we will no longer need to place a factor of $K(T)$
on the right-hand side.  

To prove \eqref{glob-free}, first observe that by dualizing the second part of \eqref{eq:compact} we have
$$ \| \int_0^{T-\tau} e^{itH_0} F_{glob}(t)\ dt \|_{\dot H^{-1/2}(M)}
\leq C (\int_0^{T-\tau} \| F_{glob}(t) \|_{H^{0,2}(M)}^2\ dt)^{1/2}$$
(for instance), so it will suffice to show that
$$ \| F_{glob}(t) \|_{H^{0,2}(M)} \leq C_{\eps_0,R,\tau} \| F(t+\tau)\|_{L^2(M)}$$
for all $0 \leq t \leq T-\tau$.  On the other hand, we know from inspection of the symbol
of $P_{med} \nabla_j \varphi$ that
$$ \|P_{med} \nabla_j \varphi F^j(t+\tau)\|_{H^{2,2}(M)} \leq C_{\eps_0} \| F(t+\tau) \|_{L^2(M)}$$
(for instance), so it will suffice by \eqref{flg-def} to show that
$$ \| e^{i\tau H} \|_{H^{2,2}(M) \to H^{0,2}(M)} \leq C_{\tau}.$$
But this can be easily established by standard energy methods\footnote{For instance, if $f \in H^{2,2}(M)$, one can first establish uniform bounds on $e^{i\tau H} f$ in $H^{2,0}(M)$ by spectral methods, then use energy methods to control $e^{i\tau H} f$ in $H^{1,1}(M)$, and then finally in $H^{0,2}(M)$, losing polynomial factors of $\tau$ in each case.  One could also argue using positive commutator methods based on \eqref{heisenberg} for such operators as $A = \langle x \rangle^4$.  We omit the details.}.  This proves \eqref{glob-free}.  Thus to prove \eqref{glob-duhamel} it suffices
to show that
\begin{equation}\label{glob-duhamel-err}
 \| \int_0^{T-\tau} (e^{itH}-e^{itH_0}) F_{glob}(t)\ dt \|_{\dot H^{-1/2}(M)}
\leq o_{R \to \infty; \eps_0}( 1 + K(T) ) (\int_0^{T-\tau} \| F(t+\tau)\|_{L^2(M)}^2\ dt)^{1/2}.
\end{equation}

At this stage it is necessary to decompose $F_{glob}(t)$
further into ``incoming'' and ``outgoing'' components, which roughly correspond to the regions of phase space where
$x \cdot \xi < 0$ and $x \cdot \xi > 0$ respectively.  Semiclassically, we expect $F_{glob}(t)$ to be supported almost
entirely in the ``incoming'' region of phase space, since it is currently far away from the origin, but came by propagating
a localized function backwards in time from $t+\tau$.  However, if this function is supported in the incoming
region of phase space, then by moving further backwards in time by $t$ it should move even further away from the origin, and in particular it should evolve much like the Euclidean flow (i.e. it should become small when $e^{itH}-e^{itH_0}$
is applied).

We now make this intuition precise.  The first step is to formalize the decomposition into incoming and outgoing waves.
We first take advantage of the fact that the spectral support of $P_{med}$ vanishes near zero.  From \eqref{flg-def} we have
\begin{equation}\label{fglob-split} 
F_{glob}(t) = (1-\varphi_R) H^2 \tilde F_{glob}(t)
\end{equation}
where
\begin{equation}\label{tfg}
 \tilde F_{glob}(t) := e^{i\tau H} \tilde P_{med} \nabla_j \varphi F^j(t+\tau)
\end{equation}
and $\tilde P_{med} := H^{-2} P_{med}$.  This factor of $H^2$ we have extracted from $P_{med}$ shall be helpful for managing the very low frequencies in the proof of Proposition \ref{psd-prop} below, 
which would otherwise cause a significant problem for this portion of the evolution
(at least in three dimensions; this step appears to be unnecessary in five and higher dimensions, for reasons similar
to why resonances do not occur in those dimensions).

We now need the following phase space decomposition associated to the \emph{Euclidean} flow $e^{-itH_0}$.

\begin{proposition}[Phase space decomposition]\label{psd-prop}  There exist operators $P_{in}, P_{out}$ such that
\begin{equation}\label{pmm} (1 - \varphi_R) H^2 = (1 - \varphi_R) H^2 P_{in} + (1 - \varphi_R) H^2 P_{out}
\end{equation}
and for which we have the estimates\footnote{The decay weights here should not be taken too seriously: indeed since we are assuming $H$ to be a compactly supported perturbation of $H_0$ we have enormous flexibility with these weights.}
\begin{equation}\label{in-propagate}
\| e^{isH_0} (1 - \varphi_R) H^2 P_{in} f \|_{H^{2,-8}(M)} 
\leq C (R^2 + |s|)^{-1-\sigma} \| f \|_{H^{20}(M)}
\end{equation}
and
\begin{equation}\label{out-propagate}
\| \langle x \rangle^{3/2+\sigma} (1 - \varphi_R) H^2 P_{out} e^{isH_0} f \|_{L^2(M)}
\leq C (R^2 + |s|)^{-1-\sigma} \| f\|_{H^{18,11}(M)}
\end{equation}
for any time $s > 0$ and all Schwartz $f$.
\end{proposition}

The proof of this proposition is a straightforward application of the principle of stationary phase, and can be justified heuristically 
by appealing to the intuition of microlocal analysis and the uncertainty principle.  It is however, a little technical 
and will be deferred to the next section.  Assuming it
for the moment, let us conclude the proof of Proposition \ref{medfreq}.  It suffices to prove \eqref{glob-duhamel}.
From \eqref{fglob-split} and Proposition \ref{medfreq} we can split $F_{glob}$ as
\begin{equation}\label{in-and-out}
 F_{glob}(t) = (1 - \varphi_R) H^2 P_{in} \tilde F_{glob}(t) + (1 - \varphi_R) H^2 P_{out} \tilde F_{glob}(t),
\end{equation}
and treat the components separately.

\divider{Term 2(a): Contribution of the global incoming part}

We now control the contribution of the incoming component of \eqref{in-and-out} to \eqref{glob-duhamel}.
We use Duhamel's formula to write
$$(e^{itH}-e^{itH_0}) (1 - \varphi_R) H^2 P_{in} \tilde F_{glob}(t)
= i \int_0^t e^{i(t-s)H} (H-H_0) e^{isH_0}
(1 - \varphi_R) H^2 P_{in} \tilde F_{glob}(t)\ ds
$$
and so it suffices to show that
\begin{align*}
 \| \int_0^{T-\tau} &\int_0^t e^{i(t-s)H} (H-H_0) e^{isH_0}
(1 - \varphi_R) H^2 P_{in} \tilde F_{glob}(t)\ ds dt \|_{\dot H^{-1/2}(M)}\\
&= o_{R \to \infty; \eps_0}( K(T) ) (\int_0^{T-\tau} \| F(t+\tau)\|_{L^2(M)}^2\ dt)^{1/2}.
\end{align*}
Substituting $t' := t-s$ using Minkowski's inequality, we can estimate
\begin{align*}
 \| \int_0^{T-\tau} &\int_0^t e^{i(t-s)H} (H-H_0) e^{isH_0}
(1 - \varphi_R) H^2 P_{in} \tilde F_{glob}(t)\ ds dt \|_{\dot H^{-1/2}(M)}\\
&\leq \int_0^{T-\tau} \| \int_0^{T-\tau-s} e^{it'H}
(H-H_0) e^{isH_0} P_{in} \tilde F_{glob}(t'+s)\ dt' \|_{\dot H^{-1/2}(M)}.
\end{align*}
Applying \eqref{Kdual}, it thus suffices to show that
\begin{equation}\label{inp}
\begin{split}
 \int_0^{T-\tau} &
\int_0^{T-\tau-s} \| \langle x \rangle^{3/2+\sigma} 
(H-H_0) e^{isH_0} (1 - \varphi_R) H^2 P_{in} \tilde F_{glob}(t'+s) \|_{L^2(M)} 
\ dt' \\
&= o_{R \to \infty; \eps_0}( (\int_0^{T-\tau} \| F(t+\tau)\|_{L^2(M)}^2\ dt)^{1/2} ).
\end{split}
\end{equation}
On the other hand, since $H-H_0$ is a compactly supported second order operator we have
$$ \| \langle x \rangle^{3/2+\sigma} 
(H-H_0) e^{isH_0} (1 - \varphi_R) H^2 P_{in} \tilde F_{glob}(t'+s) \|_{L^2(M)} 
\leq C \| e^{isH_0} (1 - \varphi_R) H^2 P_{in} \tilde F_{glob}(t'+s) \|_{H^{2,-8}(M)} $$
and so from \eqref{in-propagate} we thus have
$$
\| \langle x \rangle^{3/2+\sigma} (H-H_0) e^{isH_0} (1-\varphi_R) H^2 P_{in} \tilde F_{glob}(t'+s) \|_{L^2(M)} 
\leq C (R^2 + |s|)^{-1-\sigma} \| \tilde F_{glob}(t'+s) \|_{H^{20}(M)}.$$
From \eqref{tfg} we note that $\| \tilde F_{glob}(t'+s) \|_{L^2(M)} \leq C_{\eps_0} \| F \|_{L^2(M)}$
(noting that $\tilde P_{med}$ maps $H^{-1}(M)$ to $H^{20}(M)$ with an operator norm of $C_{\eps_0}$.  Thus
we have
\begin{align*}
\int_0^{T-\tau-s} &\| \langle x \rangle^{3/2+\sigma} 
(H-H_0) e^{isH_0} (1-\varphi_R) H^2 P_{in} \tilde F_{glob}(t'+s) \|_{L^2(M)} \\
&\leq C_{\eps_0} (R^2 + |s|)^{-1-\sigma}
(\int_0^{T-\tau} \| F(t+\tau)\|_{L^2(M)}^2\ dt)^{1/2}
\end{align*}
and the claim \eqref{inp} follows upon integrating in $s$.

\divider{Term 2(b): Contribution of the global outgoing part}

We now control the contribution of the outgoing component of \eqref{in-and-out} to \eqref{glob-duhamel}.
From \eqref{Kdual} we have
$$
  \| \int_0^{T-\tau} e^{itH} (1 - \varphi_R) H^2 P_{out} \tilde F_{glob}(t) \ dt \|_{\dot H^{-1/2}(M)}
\leq K(T) (\int_0^{T-\tau} \| S(t) \|_{L^2(M)}^2\ dt)^{1/2}
$$
where
$$ S(t) := \langle x \rangle^{3/2+\sigma} (1 - \varphi_R) H^2 P_{out} \tilde F_{glob}(t).$$
From \eqref{eq:compact}, a similar estimate holds when $e^{itH}$ is replaced by the free flow $e^{itH_0}$,
with $K(T)$ replaced by a constant $C$.  Thus to control 
this contribution to \eqref{glob-duhamel} it will suffice to show that
\begin{equation}\label{out-part}
(\int_0^{T-\tau} \| S(t) \|_{L^2(M)}^2\ dt)^{1/2}
= o_{R \to 0; \eps_0}( \int_0^{T-\tau} \| F(t+\tau) \|_{L^2(M)}^2\ dt)^{1/2}.
\end{equation}
To prove this, we first use \eqref{tfg} to expand
\begin{equation}\label{S-expand}
\begin{split} 
\| S(t) \|_{L^2(M)}^2
&= \langle \langle x \rangle^{3/2+\sigma}
(1 - \varphi_R) H^2 P_{out} e^{i\tau H} \tilde P_{med} \nabla_j \varphi F^j(t+\tau), S(t) \rangle \\
&= - \langle F^j(t+\tau),
\varphi \nabla_j P_{med} e^{-i\tau H} W(t) \rangle\\
&\leq C \| F(t+\tau) \|_{L^2(M)} \| P_{med} e^{-i\tau H} W(t) \|_{H^{1,-20}(M)} \\
&\leq C_{\eps_0} \| F(t+\tau) \|_{L^2(M)} \| e^{-i\tau H} W(t) \|_{H^{-20,-20}(M)} 
\end{split}
\end{equation}
by the support of $\varphi$ and the smoothing properties of $P_{med}$, where
\begin{equation}\label{R-def}
W(t) := P^*_{out} H^2 (1 - \varphi_R) \langle x \rangle^{3/2+\sigma} S(t).
\end{equation}
We use Duhamel's formula to write
$$
e^{-i\tau H} W(t) = e^{-i\tau H_0} W(t) - i \int_0^\tau e^{-i(\tau-s)H} (H-H_0) e^{-isH_0} W(t)\ ds,$$
and apply $H^{-20,-20}$ norms on both sides, to obtain
$$
\| e^{-i\tau H} W(t)\|_{H^{-20,-20}(M)} \leq \| e^{-i\tau H_0} W(t) \|_{H^{-20,-20}(M)}
+ \int_0^\tau \| e^{-i(\tau-s)H} (H-H_0) e^{-isH_0} W(t)\ ds \|_{H^{-20,-20}(M)}.$$
Observe that the propagator $e^{-i(\tau-s)H}$ maps $H^{-20}(M)$ to $H^{-20}(M)$ and hence to $H^{-20,-20}(M)$, uniformly
in $\tau$ and $s$.  Also, since $H-H_0$ is a compactly supported operator we see that
$H-H_0$ maps $H^{-18,-11}(M)$ to $H^{-20}(M)$.  We thus have
$$
\| e^{-i\tau H} W(t)\|_{H^1(B(0,C))} \leq C \| e^{-i\tau H_0} W(t) \|_{H^{-18,-11}(M)}
+ C \int_0^\tau \| \| e^{-is H_0} W(t) \|_{H^{-18,-11}(M)}\ ds,$$
and thus by \eqref{R-def} and the adjoint of \eqref{out-propagate}
\begin{align*}
\| e^{-i\tau H} W(t)\|_{H^1(B(0,C))} \leq &C (R^2 + \tau)^{-1-\sigma} \| S(t) \|_{L^2(M)}
+ C \int_0^\tau (R^2 + s)^{-1-\sigma} \| S(t) \|_{L^2(M)}\\
 \leq &C_{\eps_0} R^{-\sigma} \|S(t) \|_{L^2(M)}.
\end{align*}
The point here is that we have obtained a non-trivial decay in $R$.
Inserting this estimate back into \eqref{S-expand} we see that
$$ \| S(t) \|_{L^2(M)} \leq C_{\eps_0} R^{-1} \| F(t+\tau) \|_{L^2(M)},$$
and \eqref{out-part} follows.

The proof of Theorem \ref{main} is now complete, once we complete the proof of Proposition \ref{psd-prop},
which we do in the next section.

\section{Proof of the phase space decomposition}\label{psd-sec}
\def\j{{\bf j}}
\def\k{{\bf k}}
We now prove Proposition \ref{psd-prop}.  The ideas here have some similarity with a decomposition used
in \cite{st}; related ideas were also used recently in \cite{tao:focusing}.  The idea of using a decomposition
into incoming and outgoing waves to analyze the perturbation theory of the free Laplacian $H_0$ of course
goes back to Enss (and, in a different context, even earlier to the work of Lax-Phillips.)
The phase space decomposition, developed by Enss \cite{enss} and refined by Simon
\cite{simon1}, is based on the following construction. Let 
$\{\zeta_\j\}_{\j\in {\Z}^d}$ and $\{m_\j\}_{\j\in {\Z}^d}$ be 
smooth partitions of unity in ${\R}^d_x$ and ${\R}^d_\xi$ respectively with 
the property that each of the functions $\zeta_\j(x)$ and $m_\j(\xi)$ is supported
in the ball $B(\j,2)$. Define the symbols (of the $\Psi$DO's) of the projections 
$P_{in}$ and $P_{out}$
on the incoming and outgoing states: 
$$
P_{in} = \sum_{\j\cdot\k\le 0} m_\j(\xi) \, \zeta_\k(x),\qquad
P_{out}= \sum_{\j\cdot\k<0} m_\j(\xi) \, \zeta_\k(x)
$$
The following estimate, crucial to the Enss' method, reflects the expectation 
that the outgoing waves never come back to the region where they originate. 
For any
$t\ge 0, N\ge 0$ and $\j\cdot\k> 0$ such that $|\j|\ge 3$, we have
$$
\| e^{-it H_0} m_\j \zeta_\k f\|
_{L^2\big (|x|< 1/2(|\k|+t/4\big )}
\le C_N(1+t+|\k|)^{-N} \|f\|_{L^2({\R}^d) }
$$
(compare with \eqref{in-propagate} noting that the transformation 
$t\to -t$ corresponds to the change $P_{in}\to P_{out}$).
Observe that the condition $|\j|\ge 3$ ensures that we only deal with the 
outgoing waves of velocities bounded away from zero in absolute value. 

A continuous decomposition, based on coherent (Gaussian) states, was
used by Davies, \cite{davies}, while the outgoing/incoming waves defined
via projections on positive/negative spectral subspace of the dilation 
operator $x\cdot\nabla + \nabla\cdot x$ were introduced by Mourre, \cite{mourre}.

We now return to our decomposition into the incoming/outgoing wave 
needed to prove Proposition  \ref{psd-prop}.
We first observe from dyadic decomposition that it suffices
to prove the claim with $(1-\varphi_R)$ replaced by $(\varphi_{2R} - \varphi_R)$, since
the original claim then follows by replacing $R$ by $2^m R$ and summing the telescoping series
over $m \geq 0$.  By further decomposition, and some rotation and scaling, we may replace 
$\varphi_{2R} - \varphi_R$ by a smooth cutoff function $\psi = \psi_R$ supported on the ball
$B := \{ x: |x-Re_1| \leq R/100 \}$, where $e_1$ is a basis vector of $\R^n$.

Our construction will be based on Euclidean tools such as the Fourier transform, so we shall now use
the Euclidean inner product and Euclidean Lebesgue measure instead of the counterparts corresponding to
the metric $g$.  In particular the operator $H$ will no longer be self-adjoint, but this will not concern us as
we shall soon break it down into components anyway.  To reflect this change of perspective we shall write our
manifold $M$ now as $\R^n$.

We begin with the Fourier inversion formula
$$ \psi H^2 f(x) = 
\psi H^2 \tilde \psi f(x) = 
\psi H^2 \int_{\R^n} \int_{\R^n} e^{2\pi i (x-y) \cdot \xi} 
\tilde \psi (y) f(y)\ dy d\xi,$$
valid for all Schwartz $f$, where $\tilde \psi$ is a smooth cutoff to the ball $\tilde B := \{ x: |x-Re_1| \leq R/50\}$
which equals 1 on $B$.  We then split the $\xi$ integration into subspaces $\pm \xi_1 > 0$, defining
$$ P_{in} f(x) = \int_{\xi_1 < 0} \int_{\R^n} e^{2\pi i (x-y) \cdot \xi} 
\tilde \psi(y) f(y)\ dy d\xi$$
and 
$$ P_{out} f(x) = \int_{\xi_1 > 0} \int_{\R^n} e^{2\pi i (x-y) \cdot \xi} 
\tilde \psi(y) f(y)\ dy d\xi$$
where $\xi_1 := \xi \cdot e_1$ is the $e_1$ component of $\xi_1$.  We remark that these operators
are essentially Hilbert transforms in the $e_1$ direction; the multiplier is of course discontinuous in the $\xi_1$
variable, but we will never integrate by parts in this variable so this will not be a difficulty.

Clearly we have the decomposition \eqref{pmm}.  It remains to prove \eqref{out-propagate}, \eqref{in-propagate}.

\divider{Proof of \eqref{out-propagate}}

We now prove the estimate \eqref{out-propagate}.  Since $H=H_0$ on the support of $\psi$, we may write
$$ \langle x \rangle^{3/2+\sigma} \psi H^2 = R^{-5/2 + \sigma} a(x) (R \nabla_x)^4$$
for some bounded functions tensor $a(x)$ supported on $B$.  We shall think of $R\nabla_x$ as a normalized gradient on $B$.

In light of the above decomposition, it thus suffices to show that
$$ \| (R\nabla_x)^4 \int_{\xi_1 > 0} \int_{\R^3} e^{2\pi i (x-y) \cdot \xi} 
\tilde \psi(y) e^{isH_0} f(y)\ dy d\xi \|_{L^2_x(B)}
	\leq C R^{\frac{5}{2} - \sigma} (R^2 + s)^{-1-\sigma} \| f \|_{H^{18,11}(\R^3)}$$
for all all times $s > 0$.  

We first dispose of the derivatives $(R\nabla_x)^4$.
Each $x$ derivative in $R\nabla_x$ hits the phase $e^{2\pi i (x-y) \cdot \xi}$, where it
can be converted to a $y$ derivative, which after integration
by parts either hits $\tilde \psi$ or $e^{isH_0} f(y)$.  Since partial derivatives commute with
$e^{isH_0}$, we thus reduce to showing that
\begin{equation}\label{18}
 \| \int_{\xi > 0} \int_{\R^3} e^{2\pi i (x-y) \cdot \xi} 
\psi^{(j_1)}(y) e^{isH_0} f^{(j_2)}(y)\ dy d\xi \|_{L^2_x(B)}
\leq C R^{\frac{5}{2} - \sigma} (R^2 + s)^{-1-\sigma} \| f \|_{H^{18,11}(\R^3)}
\end{equation}
whenever $0 \leq j_1, j_2$ with $j_1+j_2 = 4$, where $\psi^{(j_1)} := (R \nabla_x)^{j_1} \tilde \psi$
is a minor variant of the cutoff $\tilde \psi$, and $f^{(j_2)} := (R \nabla_x)^{j_2} f$.

Applying a smooth cutoff to $f$, we can divide into two cases: the local case where $f$ is supported on the ball $\{ |x| \leq R/100\}$, or the global case where $f$ is supported on the exterior ball $\{ x \geq R/200 \}$.  Let us consider the global case first, which is rather easy and for which one can be somewhat careless with powers of $R$.  
We first observe from Plancherel (or the $L^2$ boundedness of the Hilbert transform)
that
$$
 \| \int_{\xi > 0} \int_{\R^3} e^{2\pi i (x-y) \cdot \xi} g(y)\ dy d\xi \|_{L^2_x(B)}
\leq C \| g \|_{L^2(\R^3)}.
$$
Thus to prove \eqref{18} in the global case it suffices to show that
$$ \| \psi^{(j_1)} e^{isH_0} f^{(j_2)}\|_{L^2(\R^3)} \leq C R^{\frac{5}{2} - \sigma} 
(R^2 + s)^{-1-\sigma} \| f \|_{H^{18,11}(\R^3)}.$$
We divide further into two sub-cases, the short-time case $s \leq R^2$ and the long-time case $s \geq R^2$.
In the short-time case the claim follows since
\begin{align*}
\| \psi^{(j_1)} e^{isH_0} f^{(j_2)} \|_{L^2(\R^3)} 
&\leq C \| e^{isH_0} f^{(j_2)} \|_{L^2(\R^3)} \\
&= C \| f^{(j_2)} \|_{L^2(\R^3)} \\
&\leq C R^{j_2} R^{-11} \| f \|_{H^{18,11}(\R^3)}\\
&\leq C R^{-7} \| f \|_{H^{18,11}(\R^3)}
\end{align*}
which is certainly acceptable (if $\sigma$ is small enough).  In the long time case we interpolate between the
two bounds
$$ \|  \psi^{(j_1)} e^{isH_0} g \|_{L^2(\R^3)} \leq C \| e^{isH_0} g \|_{L^2(\R^3)} \leq C \| g \|_{L^2(\R^3)}$$
and
$$ \| \psi^{(j_1)} e^{isH_0} g \|_{L^2(\R^3)} \leq C R^{n/2} \| e^{isH_0} g \|_{L^\infty(\R^3)}
\leq C R^{3/2} s^{-3/2} \| g \|_{L^1(\R^3)} \leq C R^{n/2} s^{-n/2} \| g \|_{H^{0, n/2+\sigma}(\R^3)}$$
to obtain
$$ \| \psi^{(j_1)} e^{isH_0} f^{(j_2)} \|_{L^2(\R^3)} \leq C (R/s)^{-1-\sigma}
\| f^{(j_2)} \|_{H^{0,1 + C\sigma}} \leq C (R/s)^{-1-\sigma} R^{j_2} R^{-11 + 1 + C\sigma}
\| f \|_{H^{18,11}(\R^3)}$$
which is certainly acceptable.

It remains to prove \eqref{18} in the local case, when $f$ is supported on the ball $\{ |x| \leq R/100\}$.
In this case expand the fundamental solution of $e^{isH_0}$ to write
$$ \int_{\R^3} \int_{\R^3} \chi_+(\xi) e^{2\pi i (x-y) \cdot \xi} 
\psi^{(j_1)}(y) e^{isH_0} f^{(j_2)}(y)\ dy d\xi = \int K_s(x,z) f^{(j_2)}(z)\ dz$$
where
$$ K_s(x,z) := C s^{-3/2} 
\int_{\xi_1 > 0} \int_{\R^3} e^{2\pi i (x-y) \cdot \xi} 
\psi^{(j_1)}(y) e^{-i|y-z|^2/2s}\ dy d\xi.$$
Observe from all the spatial cutoffs that $x$, $y$ are localized to the ball $\tilde B$, 
while $z$ is localized to the ball $|z| \leq R/100$.  Also, $\xi$ is localized to the half-plane
$\xi_1 > 0$.  Our task is to show that
\begin{equation}\label{remnant}
\| \int_{|x| \leq R/100} K_s(x,z) f^{(j_2)}(z)\ dz \|_{L^2(B)} \leq C R^{\frac{5}{2} - \sigma} (R^2 + s)^{-1-\sigma} \| f \|_{H^{18,11}(\R^3)}.
\end{equation}

We split $K_s$ into $K_s^{hi}$ and $K_s^{lo}$, corresponding to the regions $|\xi| \geq R^{-1+\sigma}$
and $|\xi| < R^{-1+\sigma}$ of frequency space respectively.  The contribution of $K_s^{hi}$ will be very small.
Indeed, for any $|\xi| \geq R^{-1+\sigma}$, we can evaluate the $y$ integral using stationary phase
as follows.  Observe that the $y$ gradient of the phase
$$ 2\pi (x-y) \cdot \xi - |y-z|^2 / 2s$$
is equal to
$$ - 2\pi \xi - (y-z) / s.$$
From the localizations on $y$, $z$, and $\xi$ we observe that this quantity has magnitude
at least $\geq c (|\xi| + R / s)$.  One can then do repeated integration by parts in the $y_1$ variable, gaining
an $R$ every time one differentiates the $\psi^{(j_1)}$ cutoff, to obtain the bound
$$ |\int_{\R^3} e^{2\pi i (x-y) \cdot \xi} \psi^{(j_1)}(y) e^{-i|y-z|^2/2s}\ dy|
\leq C_N R^3 (R |\xi| + R^2/s)^{-N}$$
for any $N$.  Integrating over all $|\xi| \geq R^{-1+\sigma}$, we thus see that 
$K^{hi}_s(x,z)$ is bounded by $O_N( (R^2 + s)^{-N})$ for any $N$, and so this contribution to \eqref{remnant} is easily shown to be acceptable (using crude estimates on $f^{(j_2)}$).

Now we deal with the low frequency case $|\xi| \leq R^{-1+\sigma}$.  In this case we expand $f^{(j_2)}$ and integrate by 
parts to write
\begin{align*}
|\int_{|z| \leq R/100} K^{lo}_s(x,z) &f^{(j_2)}(z)\ dz| \leq C R^{j_2} |\int_{|z| \leq R/100} \nabla_x^{j_2} K^{lo}_s(x,z) f(z)\ dz|\\
&\leq C s^{-3/2} \int_{|z| \leq R/100} \int_{|\xi| \leq R^{-1+\sigma}; \xi_1 > 0}\\
&\quad |\int_{\R^3} e^{2\pi i (x-y) \cdot \xi} \psi^{(j_1)}(y) e^{-i|y-z|^2/2s}\ dy|
R^{j_2} |\xi|^{j_2} |f(z)|\ d\xi dz.
\end{align*}
We crudely bound $R^{j_2} |\xi|^{j_2}$ by $O( R^{4\sigma} )$.  From stationary phase we also
observe the estimate
$$ |\int_{\R^3} e^{2\pi i (x-y) \cdot \xi} \psi^{(j_1)}(y) e^{-i|y-z|^2/2s}\ dy| \leq C s^{3/2}$$
while by taking absolute values everywhere we also have the estimate
$$ |\int_{\R^3} e^{2\pi i (x-y) \cdot \xi} \psi^{(j_1)}(y) e^{-i|y-z|^2/2s}\ dy| \leq C R^3$$
and hence
\begin{align*}
|\int_{|z| \leq R/100} K^{lo}_s(x,z) f^{(j_2)}(z)\ dz|
&\leq C s^{-3/2} R^{4\sigma} \int_{|z| \leq R/100} \int_{|\xi| \leq R^{-1+\sigma}}  \min(s^{n/2}, R^3)
|f(z)|\ d\xi dz \\
&\leq C R^{C\sigma} s^{-3/2} \min(s^{3/2}, R^3) R^{-3} \int_{|z| \leq R/100} |f(z)|\ dz \\
&\leq C R^{C\sigma} s^{-3/2} \min(s^{3/2}, R^3) R^{-3/2} \|f\|_{L^2(\R^3)} \\
&\leq C R^{C\sigma} R^2 (R^2 + s)^{-1-\sigma} R^{-3/2} \|f\|_{H^{18,11}(\R^3)} 
\end{align*}
which is certainly acceptable.
This completes the proof of \eqref{out-propagate}.  

\divider{Proof of \eqref{in-propagate}.}

To conclude the proof of Proposition \ref{psd-prop} it suffices to prove \eqref{in-propagate} (with $1-\varphi_R$ replaced
by $\psi$, of course).  This is in a similar spirit to the proof of \eqref{out-propagate}, although the steps will
be in a somewhat permuted order.  We begin by estimating
\begin{align*}
\| e^{isH_0} \psi H^2 P_{in} f \|_{H^{2,-8}(\R^3)} 
&\leq \sum_{k=0}^2
\| \langle x \rangle^{-8} \nabla_x^k e^{isH_0} \psi H^2 P_{in} f \|_{L^2(\R^3)} \\
&\leq \sum_{k=0}^2
\| \langle x \rangle^{-8} e^{isH_0} \nabla_x^k \psi H^2 P_{in} f \|_{L^2(\R^3)}.
\end{align*}
Now observe that
$$\nabla_x^k \psi H^2 = \sum_{j=0}^{4} R^{j - 4} a_{j,k}(x) \nabla_x^{j+k}$$
for some smooth cutoff functions $a_{j,k}$ adapted to the ball $B$.  Our task is thus to show that
$$ \| \langle x \rangle^{-8} e^{isH_0} a_{j,k} \nabla_x^{j+k} P_{in} f \|_{L^2(\R^3)}
\leq C R^{4-j} (R^2 + |s|)^{-1-\sigma} \| f \|_{H^{20}(\R^3)}$$
for all $0 \leq j \leq 4$ and $0 \leq k \leq 2$.

Fix $j,k$.  Let us first control the contribution of the weight $\langle x \rangle^{-8}$ arising
from the region $|x| \geq R/100$.  This is contribution dealt with differently in the short-time case
$s \leq R^2$ and the long-time case $s > R^2$.  In the short-time case $s \leq R^2$
we control this contribution can be controlled by
$$ \leq C R^{-8} \| e^{isH_0} a_{j,k} \nabla^{j+k} P_{in} f \|_{L^2(\R^3)}.$$
Since $e^{isH_0}$ and $a_{j,k}$ are both bounded on $L^2$, we can bound this by
$$ \leq C R^{-8} \| \nabla^{j+k} P_{in} f \|_{L^2(\R^3)}.$$
Now observe from Plancherel that $P_{in}$ is bounded on $H^{20}$ (since multiplication by $\tilde \psi$ is
certainly bounded on $H^{20}$), and so we can bound this
by $O( R^{-8} \|f\|_{H^{20}} )$, which is acceptable.  In the long-time case, we control the contribution
instead by
\begin{equation}\label{rmess}
 \leq C R^{-8+\frac{3}{2}} \| e^{isH_0} a_{j,k} \nabla^{j+k} P_{in} f \|_{H^{0,-3/2}(\R^3)}.
\end{equation}
Now we interpolate between the energy estimate
$$ \| e^{isH_0} g \|_{L^2(\R^3)} = \| g \|_{L^2(\R^3)}$$
and the decay estimate
$$ \| e^{isH_0} g \|_{H^{0,-n/2-\sigma}(\R^3)} \leq C \| e^{isH_0} g \|_{L^\infty(\R^3)}
\leq C s^{-3/2} \| g \|_{L^1(\R^3)} \leq C s^{-3/2} \| g \|_{H^{0,n/3+\sigma}(\R^3)}$$
to obtain
$$ \| e^{isH_0} g \|_{H^{0,-3/2}(\R^3)} \leq C s^{-3/2+\sigma} \| g \|_{H^{0,3/2}(\R^3)}.$$
The operator $a_{j,k}$ maps $L^2$ to $H^{0,3/2}$ with a bound of $O(R^{3/2})$, so we can therefore bound
\eqref{rmess} by
$$ \leq C R^{-8+\frac{3}{2}} s^{-3/2+\sigma} R^{3/2} \| \nabla^{j+k} P_{in} f \|_{L^2(\R^3)}.$$
As in the long-term case we can bound $ \| \nabla^{j+k} P_{in} f \|_{L^2(\R^3)}$ by $\|f\|_{H^{20}(\R^3)}$, and
so this case is also acceptable.

It remains to control the contribution in the region $|x| \leq R/100$.  We then expand out the fundamental solution
of $e^{isH_0}$, and reduce to showing that
\begin{equation}\label{fund} 
\begin{split}
s^{-3/2}
\| \langle x \rangle^{-8} \int_{\R^3} &\int_{\xi_1 < 0} \int_{\R^3}
e^{-i|x-y|^2/2s} a_{j,k}(y) \\
&\nabla_y^{j+k} e^{2\pi i (y-z) \cdot \xi} \psi(z) f(z)\ dz d\xi
dy \|_{L^2_x(|x| \leq R/100)}
\leq C R^{4-j} (R^2 + |s|)^{-1-\sigma} \| f \|_{H^{20}(\R^3)}.
\end{split}
\end{equation}
Note that the $x$ variable is localized to the ball $|x| \leq R/100$, while $y$ and $z$ are localized to the ball
$\tilde B$, and $\xi$ is localized to the half-space $\xi_1 < 0$.  Once again, we
split into the high frequencies $|\xi| \geq R^{-1+\sigma}$ and low frequencies $|\xi| < R^{-1+\sigma}$.  In
the case of the high frequencies, we move the $y$ derivatives onto $\psi(z) f(z)$ by integration by parts,
and reduce to showing that
\begin{align*}
 s^{-3/2}
\| \langle x \rangle^{-8} \int_{|\xi| \geq R^{-1+\sigma}; \xi_1 < 0} &\int_{\R^3} (\int_{\R^3}
e^{-i|x-y|^2/2s} a_{j,k}(y) e^{2\pi i (y-z) \cdot \xi}\ dy) \\
&g(z)\ dz d\xi
 \|_{L^2_x(|x| \leq R/100)}
\leq C R^{4-j} (R^2 + |s|)^{-1-\sigma} \| g \|_{L^2(\tilde B)}
\end{align*}
for some $g$.  But the $y$ derivative of the phase
$$ - |x-y|^2/2s + 2\pi (y-z) \cdot \xi$$
is
$$ - (x-y)/s + 2\pi \xi$$
which has magnitude at least $c( |\xi| + R/s)$, by the localizations on $x, y, \xi$.  Thus by
stationary phase we have
$$ \int_{\R^3}e^{-i|x-y|^2/2s} a_{j,k}(y) e^{2\pi i (y-z) \cdot \xi}\ dy = O_N( R^3 (|\xi| + R/s)^{-N})$$
for any $N > 0$, and the claim is now easy to establish by crude estimates.  Thus it remains to prove \eqref{fund}
in the low frequency case $|\xi| < R^{-1+\sigma}$.  In this case we convert the $\nabla_y^{j+k}$ derivative
to $O(|\xi|^{j+k}) = O( R^{-(j+k)} R^{C\sigma})$, and thus estimate the left-hand side of this contribution to \eqref{fund} by
\begin{align*}
\leq C s^{-3/2} &R^{-(j+k)} R^{C\sigma} R^{-3} \| \langle x \rangle^{-8} \int_{\R^3} \\
&\sup_{|\xi| \leq R^{-1+\sigma}; \xi_1 < 0}
|\int_{\R^3} e^{-i|x-y|^2/2s} a_{j,k}(y) e^{2\pi i (y-z) \cdot \xi}\ dy| \psi(z) |f(z)|\ dz \|_{L^2_x(|x| \leq R/100)}.
\end{align*}
But by stationary phase as before, we have the estimates
$$ |\int_{\R^3} e^{-i|x-y|^2/2s} a_{j,k}(y) e^{2\pi i (y-z) \cdot \xi}\ dy| \leq C \min( s^{3/2}, R^3 )$$
and thus we can bound the previous expression by
$$ \leq C s^{-3/2} R^{-(j+k)} R^{C\sigma} R^{-3} \min(s^{3/2}, R^3)
\|  \langle x \rangle^{-8} \int_{\R^3} \psi(z) |f(z)|\ dz \|_{L^2(|x| \leq R/100)}.$$
We crudely bound the $L^2(|x| \leq R/100)$ norm of $\langle x \rangle^{-8}$ by
$O( 1 )$, and use Cauchy-Schwarz we can bound the previous expression by
$$ \leq C s^{-3/2} R^{-(j+k)} R^{C\sigma} R^{-3} \min(s^{3/2}, R^3) 
 R^{3/2} \| f \|_{L^2(\R^3)},$$
which is acceptable (treating the $s \leq R^2$ and $s > R^2$ cases separately).
The proof of Proposition \ref{psd-prop} is now complete.  \endprf

\end{document}